%% file: paper.tex
\documentclass[a4paper,3p,oneside,onecolumn,preprint]{elsarticle}
\journal{Journal of Computational Physics}
\usepackage[utf8x]{inputenc}
\usepackage[english]{babel,isodate}
\usepackage{pgfplotstable}
\usepackage{graphicx}
\usepackage{grffile}
\usepackage{amsmath,amsthm}
\usepackage{amssymb}
\usepackage{todonotes}
\usepackage[integrals]{wasysym}
\usepackage{bibgerm}
\usepackage{hyperref}
\usepackage{listings}
\usepackage{environ}
\usepackage{enumerate}
\usepackage{float}
\restylefloat{table}

\begin{document}

\newcommand{\dl}{\mathrm{d}}
\newcommand{\interior}{\mathrm{int}}
\newcommand{\exterior}{\mathrm{ext}}
\newcommand{\lef}{\left\|}
\newcommand{\ri}{\right\|}
\newcommand{\diag}{\mathrm{diag}}
\newcommand\todoblock[2][]{\todo[inline,caption={},#1]{#2}}
\NewEnviron{todoenv}[1][]{\todoblock[#1]{\BODY}}
\newcommand\THsimplex[2]{$\mathcal{P}_{#1}/\mathcal{P}_{#2}$}
\newcommand\THcube[2]{$\mathcal{Q}_{#1}/\mathcal{Q}_{#2}$}
\newcommand\RT{\textup{RT}}
\newcommand\THcubeRT[3]{$\mathcal{Q}_{#1}/\mathcal{Q}_{#2}/\textup{RT}_{#3}$}
\newcommand\normal{n}

\theoremstyle{definition}
\newtheorem{theorem}{Theorem}
\newtheorem{lemma}{Lemma}
\newtheorem*{theorem*}{Theorem}
\newtheorem{observation}{Observation}
\newtheorem{remark}{Remark}
\newtheorem{algorithm}{Algorithm}


\begin{frontmatter}

\title{A Stable and High-Order Accurate Discontinuous Galerkin Based Splitting Method
for the Incompressible Navier-Stokes Equations}
\author{Marian Piatkowski\corref{cor1}}
\ead{marian.piatkowski@iwr.uni-heidelberg.de}
\author{Steffen Müthing\corref{}}
\ead{steffen.muething@iwr.uni-heidelberg.de}
\author[rvt]{Peter Bastian\corref{}}
\ead{peter.bastian@iwr.uni-heidelberg.de}
\address[rvt]{Interdisciplinary Center for Scientific Computing,
Heidelberg University, D-69120 Heidelberg}

\cortext[cor1]{Corresponding author}

\begin{abstract}
In this paper we consider discontinuous Galerkin (DG) methods
for the incompressible Navier-Stokes equations in the framework of projection methods.
In particular we employ symmetric interior penalty DG methods within the
second-order rotational incremental pressure correction scheme.
The major focus of the paper is threefold:
i) We propose a modified upwind scheme based on the Vijayasundaram numerical
flux that has favourable properties in the context of DG.
ii)
We present a novel postprocessing technique in the Helmholtz projection step
based on $H(\text{div})$ reconstruction of the pressure correction that is computed
locally, is a projection in the discrete setting and ensures that the
projected velocity satisfies the discrete continuity equation exactly. As a consequence it
also  provides local mass conservation of the projected velocity.
iii) Numerical results demonstrate the properties of the scheme for different polynomial
degrees applied to two-dimensional problems with known solution as well as
large-scale three-dimensional problems. In particular we address second-order convergence in
time of the splitting scheme as well as its long-time stability.
\end{abstract}


\end{frontmatter}


\section{Introduction}
\label{sec:introduction}

The application of discontinuous Galerkin (DG) methods to the
Navier-Stokes equations is popular due to their potentially high order of convergence,
the inf-sup stability and local mass conservation
property \cite{girault_riviere,Girault_adiscontinuous,Ferrer2011224}.
The latter is generally not fulfilled for conforming
finite element discretizations.
In addition to the $2 \times 2$ block structure arising from the saddle point system
discontinuous Galerkin methods offer a further block
structure when the unknowns associated with one cell of the mesh are grouped together.
This data structure is essential for high-performance implementations
of the discontinuous Galerkin method \cite{2016arXiv160701323K,sum-fact-dg}
as it avoids costly memory gather and scatter operations when compared to
conforming finite element methods.

Operator splitting methods for solving the instationary Navier-Stokes
equations has been subject to detailed investigations for the recent
decades. One possibility in the splitting methods is to split between
the convective term and the saddle point structure which is realized in
Glowinski's $\Theta$-scheme, \cite{CBO9780511574856A009,Elman}.
Another possibility is to split between
incompressibility and dynamics which has been independently developed by
Chorin \cite{ChorinProjectionNavierStokes} and T\'{e}mam
\cite{TemamProjectionNavierStokes} and is referred to as Chorin's
projection method. The latter splitting schemes have the appealing
feature that at each time step, instead of solving a saddle point
system, one only has to solve a vector-valued heat equation for the
velocity (in the Stokes case) and a Poisson equation for the
pressure. The choice of artificial boundary conditions on the
pressure Poisson equation is a delicate issue in
projection methods of this class \cite{Rannacher1992,WeinanE1,WeinanE2}.
Several higher-order extensions of Chorin's first order method have been suggested in
the literature
\cite{KARNIADAKIS1991414,FLD:FLD373,GuermondShen,Guermond04onthe,Karniadakis2005,Guermond06}.
Here we concentrate on the classic incremental pressure-correction scheme (IPCS) \cite{GODA197976}
and the rotational incremental pressure-correction scheme (RIPCS) \cite{FLD:FLD373}.

The use of a DG spatial discretization within splitting schemes is a current subject
of active research. A naive computation of the divergence free velocity by subtraction of
the rotation free part is reported to be unstable when the spatial mesh is
coarse and the time step is small, see \cite{Steinmoeller2013480,Joshi2016120,2016arXiv160701323K},
where several local postprocessing techniques are discussed to overcome this difficulty.
In this paper we propose a new postprocessing
technique based on $H(\text{div})$ reconstruction of the discrete pressure gradient which is popular in
porous media flow computations \cite{BastianRiviere,ErnHdiv2007}.
The new approach provides a discrete velocity that satisfies the discrete continuity equation exactly and
in consequence is locally mass conservative and defines a projection. These properties are not
satisfied by the postprocessing schemes available in the literature.

The structure of the paper is organized as follows: In section
\ref{sec:dg-discr-navi} we recapitulate the discontinuous Galerkin
discretization by the interior penalty method as presented in
\cite{Girault_adiscontinuous,girault_riviere}. In section
\ref{sec:projection-methods} we discuss the Helmholtz decomposition, prove
our main result and present the projection methods
In section \ref{sec:numer-exper} we elaborate on
numerical experiments for the discontinuous Galerkin discretization
based on the reference problems by
\cite{GuermondShen,Guermond04onthe,doi:10.1137/040604418,Guermond06} and assess
the properties of the new postprocessing scheme.

\section{Discontinuous Galerkin discretization of the incompressible Navier-Stokes equations}
\label{sec:dg-discr-navi}

In this section we present the spatial discretization of the Navier-Stokes system with
an interior penalty DG method taken from \cite{Girault_adiscontinuous}. The convective term
is discretized using the Vijayasundaram flux.

The instationary incompressible Navier-Stokes equations in
an open and bounded domain $\Omega\subset\mathbb{R}^d$ ($d=2,3$)
determining the velocity $v$ and pressure $p$ for a
right-hand side $f$, constant viscosity $\mu$ and density $\rho$ are
given by
\begin{subequations}
\label{eq:incomp_navierstokes}
\begin{align}
\rho\partial_t v - \mu\Delta v + \rho (v\cdot\nabla) v + \nabla p &= f
  && \textup{in} \; \Omega \times (0,T] \\
\nabla\cdot v &= 0 && \textup{in} \; \Omega \times (0,T]\\
  v &= v_0 && \textup{for} \; t=0 \\
\intertext{and either Dirichlet boundary condition for the velocity:}
      v &= g  && \textup{on} \; \Gamma_D = \partial\Omega, t\in (0,T] \label{eq:Dirichlet1}\\
\intertext{together with}
      \int_\Omega p \dl x &= 0 && \textup{for all} \; t\in (0,T] \label{eq:Dirichlet2}\\
\intertext{or mixed boundary conditions:}
      v &= g && \textup{on} \; \Gamma_D\neq\{\partial\Omega, \emptyset\} \label{eq:Mixed1}\\
      \mu\nabla v \normal - p\normal &= 0 &&\textup{on} \; \Gamma_N = \partial\Omega\setminus \Gamma_D \label{eq:Mixed2}
\end{align}
\end{subequations}
with $(0,T]$ being the time interval of interest. For pure Dirichlet boundary conditions $g$ is required
to satisfy the compatibility condition $\int_{\partial\Omega} g\cdot \normal \dl x=0$.
In the numerical examples below we will also
consider periodic boundary conditions in addition.
Under appropriate assumptions the Navier-Stokes problem in weak form
has a solution $(v,p)$ in $(H^1(\Omega))^d\times L^2(\Omega)$ for $t\in(0,T]$, \cite{GiraultRaviartBook,TemamBook}.
In case of pure Dirichlet boundary conditions the pressure is only
determined up to a constant and is in the space $L^2_0(\Omega) = \{q\in L^2(\Omega) \mid \int_\Omega q \dl x = 0  \}$.

For the discretization let $\mathcal{E}_h$ be an affine cubic mesh
(the restriction to affine meshes is only needed when the Raviart-Thomas reconstruction
is used) with maximum diameter $h$. We denote by $\Gamma_h^\interior$
the set of all interior faces, by $\Gamma_h^{D}$ the set of all faces intersecting
with the Dirichlet boundary $\Gamma_D$ and by $\Gamma_h^N$ the set of all faces intersecting with the
mixed boundary $\Gamma_N$. We set $\Gamma_h = \Gamma_h^\interior
\cup \Gamma_h^D \cup \Gamma_h^N$. To an interior face $e\in\Gamma_h^\interior$
shared by elements $E_e^1$ and $E_e^2$ we define an orientation through
its unit normal vector $\normal_e$ pointing from $E_e^1$ to $E_e^2$. The
jump and average of a scalar-valued function $\phi$ on a face
is then defined by
\begin{align}
  [\phi] &= \phi\mid_{E_e^1} - \phi\mid_{E_e^2} \quad \; = \phi^\interior -
  \phi^\exterior, \label{eq:jump_average} \\
  \{\phi\} &= \frac12 \phi\mid_{E_e^1} + \frac12 \phi\mid_{E_e^2} =
  \frac12 \phi^\interior + \frac12 \phi^\exterior \notag \; .
\end{align}
Note that the definition of jump and average can be extended in a
natural way to vector and matrix-valued functions.
If $e\in\partial\Omega$ then $\normal_e$ corresponds to the outer normal
vector $\normal$.
Below we make heavy use of the identities and notation, respectively:
\begin{align}
  [u v] &= [u] \{v\} + \{u\} [v] \;,
  & (u,v)_{0,\omega} &= \int_{\omega} u v\,dx\; ,
  &&\text{($u,v$ scalar-valued)} \label{eq:jump_of_product}\\
  [u \cdot v] &= [u] \cdot \{v\} + \{u\} \cdot [v] \;,
  & (u,v)_{0,\omega} &= \int_{\omega} u \cdot v\,dx\; ,
  &&\text{($u,v$ vector-valued)} \notag \\
  [u : v] &= [u] : \{v\} + \{u\} : [v] \;,
  & (u,v)_{0,\omega} &= \int_{\omega} u : v\,dx\; ,
  &&\text{($u,v$ matrix-valued)} \notag .
\end{align}
The DG discretization on cuboid meshes is based on the non-conforming finite
element space of polynomial degree $p$
\begin{align}
Q_h^p = \{ v\in L^2(\Omega) \mid v|_E= q\circ\mu_E^{-1},  q\in \mathbb{Q}_{p,d}, E\in\mathcal{E}_h \}
\end{align}
where $\mu_E : \hat E \to E$ is the transformation from the reference cube $\hat E$ to $E$
and $\mathbb{Q}_{p,d}$ is the set of polynomials of maximum degree $p$ in $d$ variables.
The approximation spaces for velocity and pressure are then
\begin{subequations}
\begin{align}
X^p_h \times M_h^{p-1} &= (Q_h^p)^d \times (Q_h^{p-1}\cap L^2_0(\Omega))
&&\text{(Dirichlet b. c.)},\\
X^p_h \times M_h^{p-1} &= (Q_h^p)^d \times Q_h^{p-1}
&&\text{(mixed b. c.)} .
\end{align}
\end{subequations}
We make use of the following mesh-dependent forms
defined on $X_h^p \times X_h^p$, $X_h^p \times M_h^{p-1}$ and $M_h^{p-1}$,
respectively:
\begin{subequations}
\begin{align}
  a(u,v) &= d(u,v) + J_0(u,v), \ \text{where}\\
  d(u,v) &= \sum_{E\in\mathcal{E}_h} (\nabla u, \nabla v)_{0,E}
  -\sum_{e\in\Gamma_h^\interior} (\{\nabla u\}\normal_e, [v])_{0,e}
   - \sum_{e\in\Gamma_h^D} (\nabla u^\interior\normal_e , v^\interior)_{0,e}, \label{eq:dg_discr_diffusion}\\
  J_0(u,v) &= \epsilon \sum_{e\in\Gamma_h^\interior} (\{\nabla v\}\normal_e , [u])_{0,e}
  + \epsilon \sum_{e\in\Gamma_h^D} (\nabla v^\interior\normal_e,u^\interior)_{0,e} \notag\\
  &\quad + \sum_{e\in\Gamma_h^\interior} \frac{\sigma}{h_e} ( [u] ,
  [v])_{0,e} + \sum_{e\in\Gamma_h^D} \frac{\sigma}{h_e}
  (u^\interior , v^\interior)_{0,e}, \label{eq:dg_penalty} \\
  b(v,q) &= -\sum_{E\in\mathcal{E}_h} (\nabla\cdot v,q)_{0,E} +
  \sum_{e\in\Gamma_h^\interior} ([v]\cdot\normal_e,\{q\} )_{0,e}
  +\sum_{e\in\Gamma_h^D} (v^\interior \cdot \normal,q^\interior )_{0,e} ,\\
  l(v;t) &=\sum_{E\in\mathcal{E}_h} (f(t) , v)_{0,E}  +
  \epsilon \sum_{e\in\Gamma_h^D} (\nabla v^\interior\normal_e , g(t) )_{0,e}
  + \sum_{e\in\Gamma_D} \frac{\sigma}{h_e} (g(t),v^\interior)_{0,e} ,\\
  r(q;t) &=\sum_{e\in\Gamma_h^D}  (g(t)\cdot\normal,q^\interior)_{0,e}
  \; .
\end{align}
\end{subequations}
Here we made the time dependence of the right hand side functionals explicit. For ease of writing
this will be omitted mostly below.
In the interior penalty parameter
\(\sigma / h_e\), the denominator accounts for the mesh
dependence. The formula for \(h_e\),
\begin{equation*}
  h_e = \begin{cases}
    \frac{\min\left( \left| E^{\interior}(e) \right|, \left|
          E^{\exterior}(e) \right| \right)}
    { \left| e \right| } &, E^{\interior}(e) \cap E^{\exterior}(e) = e \\
    \frac{ \left| E^{\interior}(e) \right| }{ \left| e \right| } &,
    E^{\interior}(e)
    \cap \Gamma_D = e
  \end{cases} \; ,
\end{equation*}
has been stated in \cite{HoustonHartmann2008} where it was proven that
this choice ensures coercivity of the bilinear form for anisotropic
meshes.
For $\sigma$ we choose \(\sigma = \alpha p(p+d-1)\) as in
\cite{amg4dg} with \(\alpha\) a user-defined parameter.
In $J_0$ the \emph{Symmetric Interior Penalty Galerkin (SIPG)}
($\epsilon = -1$) method is
preferred since the matrix of the linear system in absence of the
convection term is then symmetric. Other choices are
the NIPG ($\epsilon = 1$) or IIPG ($\epsilon = 0$) method.

A first discretization of the nonlinear term in the Navier-Stokes
equations is the standard (or centered) discretization,
\begin{align}
  c(u;z,\theta) = \sum_{E\in\mathcal{E}_h} ( (u\cdot\nabla) z ,\theta)_{0,E} ,
\end{align}
that lets us define the discrete in space, continuous in time
formulation of the Navier-Stokes problem (\ref{eq:incomp_navierstokes}).
Find $v_h(t) : (0,T]\to X_h^p$, $p_h(t) : (0,T]\to M_h^{p-1}$:
\begin{subequations}
\begin{align}
  \rho(\partial_t v_h, \varphi)_{0,\Omega} + \mu a(v_h,\varphi) +
  \rho c(v_h;v_h,\varphi) + b(\varphi,p_h) &=
  l(\varphi;t) , \label{eq:variational_momentum }\\
  b(v_h,q) &= r(q;t) , \label{eq:variational_incomp}
\end{align}
\end{subequations}
for all $(\varphi,q)\in X_h^p\times M_h^{p-1}$.
This formulation for the variational form \(c\) is only applicable for
small Reynolds numbers. Therefore we present for higher Reynolds numbers
an upwind discretization in Section \ref{sec:discr-conv-part}.
The following observation will be used in several circumstances below.
\begin{remark}\label{rem:alternativebh}
The bilinear form $b(v,q)$ has the equivalent representation
\begin{equation}
b(v,q) = \sum_{E\in\mathcal{E}_h} (v,\nabla q)_{0,E}
- \sum_{e\in\Gamma_h^{\text{int}} } (\{v \}\cdot \normal_e,[q])_{0,e}
- \sum_{e\in\Gamma_h^{N}} (v\cdot \normal,q)_{0,e} .
\end{equation}
This holds true for Dirichlet and mixed boundary conditions (in the former case
just set $\Gamma_h^N=\emptyset$).
\begin{proof} Follows from integration by parts and
  \eqref{eq:jump_of_product}.
\end{proof}
\end{remark}
As a corollary we obtain the following local mass conservation property
by testing \eqref{eq:variational_incomp}
with $q=\chi_E$, the characteristic function of element $E$, and using
Remark \ref{rem:alternativebh}:
\begin{equation}
\sum_{e\in\Gamma_h^{\text{int}}\cap\partial E} (\{v \}\cdot n_e,1)_{0,e}
+ \sum_{e\in\Gamma_h^{N}\cap\partial E} (v\cdot \normal_e,1)_{0,e}
+ \sum_{e\in\Gamma_h^D\cap\partial E}  (g\cdot\normal,1)_{0,e} = 0 .
\label{eq:local_mass_conservation}
\end{equation}

\subsection{Upwind discretization of the convective part}
\label{sec:discr-conv-part}

For higher Reynolds numbers we employ a suitable upwind discretization
based on the Vijayasundaram numerical flux adapted from DG methods for
inviscid compressible flow \cite{feistauer_book,Dole2004727}.

Note that due to $\nabla\cdot v=0$ the convective term in the momentum equations can
be written equivalently as $(v\cdot\nabla) v = \nabla\cdot (v \otimes v)$
where
\begin{equation*}
F(v) = v\otimes v = [v_1 v, \ldots, v_d v] = [F_1(v), \ldots, F_d(v)]
\end{equation*}
is the convective flux matrix with columns $F_k(v)=v_kv$ and $\nabla_v F_k(v)=(v_k I + v\otimes e_k)$.
$I$ denotes the identity matrix and $(e_k)_i=\delta_{ik}$ are the coordinate unit vectors.
In order to derive the upwinding we consider the first order system
$$\partial_t v + \nabla\cdot F(v) = 0$$
which is said to be hyperbolic
if the matrix
$$ P(v,n) = \sum_{k=1}^d n_k \nabla_v F_k(v) = (v\cdot\normal)I + v\otimes n $$
is real diagonalizable for all $v, n\in\mathbb{R}^d$ with $\|\normal\|=1$ \cite{Evans}.
This is indeed the case for $v\cdot n\neq 0$. When $v\cdot n=0$, $P(v,n)=v\otimes n$ has
$d$ eigenvalues zero with a corresponding eigenspace $W_n^\perp = \{w : w\cdot n=0\}$ of dimension
$d-1$.

When discretizing the conservative form of the convective terms
with DG one uses element-wise integration by parts to arrive at
\begin{equation*}
\begin{split}
c(v;v,\varphi) &=  (\nabla\cdot F(v),\varphi)_{0,\Omega} =
  \sum_{E\in\mathcal{E}_h} (\nabla\cdot F(v) ,\varphi)_{0,E}
  = -\sum_{E\in\mathcal{E}_h} ( F(v) , \nabla\varphi)_{0,E} +
  \sum_{E\in\mathcal{E}_h} (F(v)\normal,\varphi)_{0,\partial E} \\
  &= -\sum_{E\in\mathcal{E}_h} ( F(v) , \nabla\varphi)_{0,E} +
  \sum_{e\in\Gamma_h^{\text{int}}} ([ F(v)\normal_e \cdot \varphi ],1)_{0,e}
  + \sum_{e\in\Gamma_h^{D}\cup\Gamma_h^{N}}  (F(v)\normal,\varphi)_{0,e}
  \label{eq:part_integration_convection_conservative_1}
\end{split}
\end{equation*}
Now the flux in face normal direction $F(v)\normal_e$ needs to be replaced by a \textit{consistent} and \textit{conservative}
numerical flux function $\hat F(v,\normal_e)$ which we now derive.
Since $F_k(v)=v_kv$ is homogeneous of degree 2
(i.e. $F_k(\alpha v)=\alpha^2 F_k(v)$ for $\alpha$ a real number)
it admits a representation
\begin{equation*}
\label{eq:the_A_matrix}
F_k(v) = \frac12 \nabla_v F_k(v) v
\end{equation*}
and therefore
\begin{equation*}
F(v)\normal = \frac12 P(v,n) v = \frac12 [(v\cdot\normal)I + v\otimes n] v =: B_{\frac12}(v) v\; .
\end{equation*}
Using the identity $(v\cdot n) v = (v\otimes n) v$ we see
\begin{equation*}
F(v)\normal = B_\beta(v) v := [(1-\beta) (v\cdot\normal) + \beta v\otimes n] v
\end{equation*}
for any \(\beta\in[0,1]\). For $v\cdot n\neq 0$, $B_\beta(v,\normal)$ is real diagonalizable with
eigenvalues $\lambda_{\beta,i}\in\mathbb{R}$
and a full set of right eigenvectors $r_i$, $\text{span}\{r_1,\ldots,r_{d-1}\}=W_n^\perp$, $r_d=v$,
admitting the decomposition
\begin{equation*}
B_\beta(v,\normal) = B_\beta^+(v,\normal) + B_\beta^-(v,\normal),
\end{equation*}
where $B_\beta^{\pm}(v,\normal) = T D_\beta^{\pm} T^{-1}$,
$T=[r_1,\ldots,r_d]$,
$D_\beta^\pm$ are diagonal matrices with
$(D_\beta^+)_{ii} = \max(0,\lambda_{\beta,i})$ and $(D_\beta^-)_{ii} = \min(0,\lambda_{\beta,i})$ (all eigenvectors and eigenvalues
depending on $v$ and $\normal$).

Following \cite{Dole2004727}, in the DG scheme we employ the \textit{Vijayasundaram} numerical flux given by
\begin{equation}
\label{eq:numerical_flux}
\hat F_\beta (v,n_e) = B^+_\beta(\{v\},n_e) v^\text{int} + B^-_\beta(\{v\},n_e) v^\text{ext} \; .
\end{equation}
Here the matrices $B_\beta^\pm(\{v\},n_e)$ are \textit{not} applied to $\{v\}$ and therefore
$(\{v\}\cdot n)I$ and $\{v\}\otimes n$ act differently. The effect is shown by the following
\vspace{-0.5em} 
\begin{observation} Assume $\{v\}\cdot n\neq 0$. Then the numerical flux \eqref{eq:numerical_flux}
satisfies
\begin{equation*}
\begin{split}
\hat F_\beta (v,n_e) &= (1-\beta) \left[ \max(0,\{v\}\cdot n_e) v^\text{int} + \min(0,\{v\}\cdot n_e) v^\text{ext} \right]\\
&\quad + \beta [H(\{v\}\cdot n_e)(v^\text{int}\cdot n_e) + H(-\{v\}\cdot n_e)(v^\text{ext}\cdot n_e)] \{v\} \; .
\end{split}
\end{equation*}
where $H(x)$ is the Heaviside function.
\begin{proof} We consider the interior part. The eigenvectors of $B_\beta(\{v\},n_e)$ are $d-1$ vectors
spanning $W_n^\perp$ and $\{v\}$ independent of $\beta\in[0,1]$. We can uniquely decompose
$$
v^\text{int}
 = \left( v^\text{int} - \frac{v^\text{int}\cdot n_e}{\{v\}\cdot n_e} \{v\}\right)
 + \frac{v^\text{int}\cdot n_e}{\{v\}\cdot n_e} \{v\}
 = w + \alpha \{v\}
$$
where $w\in W_n^\perp$. Now
\begin{equation*}
\begin{split}
B^+_\beta&(\{v\},n_e) v^\text{int} = B^+_\beta(\{v\},n_e)  ( w + \alpha \{v\}) \\
&= (1-\beta) \max(0,\{v\}\cdot n_e) w + \max(0,\{v\}\cdot n_e) \alpha \{v\}\\
&= (1-\beta) \max(0,\{v\}\cdot n_e) (v^\text{int} - \alpha \{v\}) + \max(0,\{v\}\cdot n_e) \alpha \{v\}\\
&= (1-\beta) \max(0,\{v\}\cdot n_e) v^\text{int}
+ \beta \frac{\max(0,\{v\}\cdot n_e) (v^\text{int}\cdot n_e)}{\{v\}\cdot n_e} \{v\} \; .
\end{split}
\end{equation*}
$B^-_\beta(\{v\},n_e) v^\text{ext}$ can be treated in the same way.
\end{proof}
\end{observation}
The observation shows that for $\beta>0$ the $v\otimes n_e$ part gives a contribution in the
flux in the direction of $\{v\}$, i.e. a central flux which moreover might have the wrong sign
since the signs of $\{v\}\cdot n_e$ and $v^\text{int}\cdot n_e$ or $v^\text{ext}\cdot n_e$
might differ since the DG velocity is not in $H(\text{div};\Omega)$. (Note, however, that the new projection
scheme to be described below improves significantly on this point).
Also note that the upwind decision is based on the average velocity which is locally mass
conservative due to \eqref{eq:local_mass_conservation}.

For these reasons we propose to employ $\beta=0$ in the numerical flux function, leading to
the simple form:
\begin{equation*}
\hat F_e(v,\normal_e) = \left\{\begin{array}{ll}
\max(0,\{v\}\cdot\normal_e) v^{\text{int}} + \min(0,\{v\}\cdot\normal_e) v^\text{ext}  & e\in\Gamma_h^\text{int}\\
\max(0,v^\text{int}\cdot\normal_e) v^{\text{int}} + \min(0,v^\text{int}\cdot\normal_e) g & e\in\Gamma_h^D\\
\max(0,v^\text{int}\cdot\normal_e) v^{\text{int}} & e\in\Gamma_h^N
\end{array}\right. \;
  \label{eq:vijayasundaram_flux}
\end{equation*}
and the upwind DG discretization of the convective term
\begin{equation}
\begin{split}
\hat{c}(v;v,\varphi) &= -\sum_{E\in\mathcal{E}_h} ( F(v) , \nabla\varphi)_{0,E} +
  \sum_{e\in\Gamma_h^{\text{int}}} (\hat{F}_e(v,\normal_e) , [\varphi ])_{0,e}
 + \sum_{e\in\Gamma_h^{D}\cup\Gamma_h^{N}}
 (\hat{F}_e(v,\normal_e),\varphi)_{0,e} \; .
  \label{eq:part_integration_convection_conservative_final}
\end{split}
\end{equation}
In the following computations we will use this variational form
\(\hat{c}\) in solving equation (\ref{eq:variational_momentum }).

\section{Projection methods}
\label{sec:projection-methods}

\subsection{Continuous Helmholtz decomposition}
\label{sec:helmh-decomp}

The Helmholtz decomposition takes a fundamental role in the construction
of splitting methods for incompressible flows. It states that any vector
field in $L^2(\Omega)^d$ can be decomposed into a divergence-free
contribution and an irrotational contribution,
see e.g. \cite{ChorinProjectionNavierStokes, Karniadakis2005,
book:941322,Schweizer2013,Bhatia:2013:HDS:2498747.2498999}.
In order to define the decomposition boundary conditions on the
pressure need to be enforced which are not part of the underlying Navier-Stokes equations.
The choice and consequence of these boundary conditions is a delicate issue in
projection methods \cite{Rannacher1992,WeinanE1,WeinanE2}.
Before turning to the Helmholtz decomposition in the discrete setting of
DG methods we recall the Helmholtz decomposition in the weak continuous setting.

First consider Dirichlet boundary conditions \eqref{eq:Dirichlet1}, \eqref{eq:Dirichlet2}. Let us
denote the space of weakly divergence free functions by
\begin{equation}
H(\Omega) := \{ v \in L^2(\Omega)^d \mid (v,\nabla f)_{0,\Omega} - (g\cdot\normal,f)_{0,\Gamma_D}
= 0 \; \forall f\in H^1(\Omega) \}
\end{equation}
where $\Gamma_D=\partial\Omega$. This definition is motivated by the identity
$(\nabla\cdot v,f)_{0,\Omega} = - (v,\nabla f)_{0,\Omega} + (g\cdot\normal,f)_{0,\Gamma_D} = 0$
which holds true for $v\in H(\text{div};\Omega) = \{u\in L^2(\Omega)^d \mid \nabla\cdot u\in L^2(\Omega)\}$.
In that case the normal component of $v$ can be prescribed on the boundary.
In addition, we employ the pressure space
\begin{equation}
\Psi_D(\Omega) := \{ q\in H^1(\Omega) \mid (q,1)_{0,\Omega} = 0 \} .
\end{equation}
in the following decomposition.

\begin{theorem}[Helmholtz decomposition, Dirichlet boundary conditions]
  \label{theo:helmholtz_decomposition}
For any $w\in L^2(\Omega)^d$ there are unique functions
$v\in H(\Omega)$ and $\psi\in \Psi_D(\Omega)$ such that
$$ w = v + \nabla \psi.$$
\begin{proof} Define $\psi\in\Psi_D(\Omega)$ by
\begin{equation}
\label{eq:Poisson1}
(\nabla\psi,\nabla q)_{0,\Omega} = (w,\nabla q)_{0,\Omega} - (g\cdot\normal,q)_{0,\Gamma_D}
\qquad \forall q\in \Psi_D(\Omega).
\end{equation}
According to the Lax-Milgram theorem this problem has a unique solution.
Since any $f\in H^1(\Omega)$ can be written as
$f=q+c$ with $q\in \Psi_D(\Omega)$
and $c$ a constant function,
equation \eqref{eq:Poisson1} holds also true for all test functions in $H^1(\Omega)$ (Note
the compatibility condition on $g$).
Now set $v=w-\nabla\psi$ and verify that
$(v,\nabla f)_{0,\Omega} - (g\cdot\normal,f)_{0,\Gamma_D} =0$ for all $f\in H^1(\Omega)$.
\end{proof}
\end{theorem}

\begin{remark}
\begin{enumerate}[1)]
\item Note that equation \eqref{eq:Poisson1} is the weak formulation of a Poisson equation with
homogeneous Neumann boundary conditions.
\item The map $\mathcal{P} : L^2(\Omega)^d \to H(\Omega)$ given by $\mathcal{P} w = w-\nabla\psi$ is a projection
since the right hand side of \eqref{eq:Poisson1} is zero for $w\in H(\Omega)$.
$\mathcal{P}$ is called the \textit{continuous Helmholtz projection}.
\item The construction above can be equivalently written as
\begin{subequations}
\label{eq:equivalent_form_decomposition}
\begin{align}
(v,\varphi)_{0,\Omega} + (\nabla\psi,\varphi)_{0,\Omega} &= (w,\varphi)_{0,\Omega} &&\forall \varphi\in L^2(\Omega)^d\\
(v,\nabla q)_{0,\Omega} &= (g\cdot\normal,q)_{0,\Gamma_D} &&\forall q\in\Psi_D(\Omega)
\end{align}
\end{subequations}
since from the first equation we get $v=w-\nabla\psi$ and inserting in the
second equation yields \eqref{eq:Poisson1}.
\item In Chorin's classical projection scheme \cite{ChorinProjectionNavierStokes} the (divergence-free)
velocity $v^{k+1}$ and pressure $p^{k+1}$ at time $t^{k+1}$ are computed from a tentative velocity $w^{k+1}$ by the system
\begin{align*}
\frac{v^{k+1} - w^{k+1}}{\Delta t} + \nabla p^{k+1} &= 0\\
\nabla\cdot v^{k+1} &= 0
\end{align*}
in strong form. Setting $\psi^{k+1}=\Delta t p^{k+1}$ this is equivalent to
\begin{align*}
v^{k+1}  + \nabla \psi^{k+1} &= w^{k+1}\\
\nabla\cdot v^{k+1} &= 0
\end{align*}
which is the strong form of \eqref{eq:equivalent_form_decomposition}. Thus, $\psi/\Delta t$ from
the Helmholtz decomposition is the new pressure from Chorin's projection scheme.
\hfill$\square$
\end{enumerate}
\end{remark}

\noindent In the case of mixed boundary conditions \eqref{eq:Mixed1}, \eqref{eq:Mixed2} the space $\Psi_D(\Omega)$ is
replaced by
\begin{equation}
\Psi_M(\Omega) := \{ q\in H^1(\Omega) \mid q=0 \text{ a.e. on $\Gamma_N$} \}
\end{equation}
employing homogeneous Dirichlet boundary conditions on $\Gamma_N$. This can
be understood from \eqref{eq:Mixed2} which implies $p\approx 0$ for small $\mu$, i.e.
large Reynolds number. The irrotational part is defined as in \eqref{eq:Poisson1} with
$\Psi_D(\Omega)$ replaced by $\Psi_M(\Omega)$, meaning that $\psi$ satisfies homogeneous
Neumann conditions on $\Gamma_D$ and homogeneous Dirichlet conditions on $\Gamma_N$.
 Again, $v\in H(\Omega)$ is uniquely
defined (observe that now $\Gamma_D\subset\partial\Omega$ in $H(\Omega)$).


\subsection{Discrete Helmholtz decomposition}
\label{sec:discr-helmh-decomp}

We now seek discrete versions $\mathcal{P}_h : X_h^p\to X_h^p$ of the
Helmholtz projection operator $\mathcal{P}$.
A direct reconstruction of the weakly divergence free velocity as $v=w-\nabla \psi$ in DG splitting schemes
is reported to be unstable when the spatial mesh is
coarse and the time step is small \cite{Steinmoeller2013480,Joshi2016120,2016arXiv160701323K} and several
local postprocessing techniques are discussed in the literature.
Here we propose a new postprocessing
technique based on $H(\text{div})$ reconstruction which is popular in porous media flows
\cite{BastianRiviere,ErnHdiv2007}. These reconstructions are element-local, easy to compute and
provide a locally mass conservative projected velocity, a property not shared by the reconstructions
in \cite{Steinmoeller2013480,2016arXiv160701323K}. \cite{Joshi2016120} takes into
account inter-element continuity in a regularized least-squares sense but does not provide a projection.
The construction presented here
is easier to compute, provides exact local mass conservation, satisfies the discrete continuity equation exactly
and provides a projection.

\subsubsection{Standard projection}

For any given tentative velocity $w_h\in X_h^p$ the straightforward translation of the
Helmholtz decomposition \eqref{eq:equivalent_form_decomposition} in the
DG setting reads
\begin{subequations}
\begin{align}
(v_h,\varphi)_{0,\Omega} + (\nabla_h\psi_h,\varphi)_{0,\Omega} &= (w_h,\varphi)_{0,\Omega}
&&\forall \varphi\in X_h^p, \label{eq:discreteHD1}\\
b(v_h,q) &= r(q) &&\forall q\in M_h^{p-1} . \label{eq:discreteHD2}
\end{align}
\end{subequations}
Note that the second equation requires the projected velocity to satisfy the discrete form
of the continuity equation \eqref{eq:variational_incomp} at fixed time
(hence silently dropping the time dependence from now).
From the first condition \eqref{eq:discreteHD1} we get $v_h + \nabla_h\psi_h=w_h \Leftrightarrow v_h = w_h - \nabla_h\psi_h$
since all involved functions are in $X_h^p$.
Inserting this into \eqref{eq:discreteHD2}
yields an equation for $\psi_h$:
\begin{equation*}
b(\nabla\psi_h,q)= b(w_h,q)-r(q) \quad \forall q\in M_{h}^{p-1} .
\end{equation*}
Using Remark \ref{rem:alternativebh} on the left hand side we get
\begin{equation}
b(\nabla\psi_h,q) =
\sum_{E\in\mathcal{E}_h} (\nabla\psi_h,\nabla q)_{0,E}
- \sum_{e\in\Gamma_h^{\text{int}} } (\{\nabla\psi_h \}\cdot \normal_e,[q])_{0,e}
- \sum_{e\in\Gamma_h^{N}} (\nabla\psi_h\cdot \normal_e,q)_{0,e} .
\label{eq:pressure_poisson2}
\end{equation}
This is part of the standard SIPG formulation of Poisson's equation with homogeneous Neumann
boundary conditions on $\Gamma^D$
with the stabilization terms missing. In order to stabilize, we define
\begin{equation}
\begin{split}
j_0(\psi_h,q) &=
- \sum_{e\in\Gamma_h^{\text{int}} } (\{\nabla q \}\cdot \normal_e,[\psi_h])_{0,e}
+ \sum_{e\in\Gamma_h^{\text{int}} } \frac{\sigma}{h_e} ([q],[\psi_h])_{0,e}
 - \sum_{e\in\Gamma_h^{N}} (\nabla q\cdot \normal_e,\psi_h)_{0,e}
+ \sum_{e\in\Gamma_h^{N}} \frac{\sigma}{h_e} (q,\psi_h)_{0,e} .
\end{split}
\end{equation}
and solve the \textit{stabilized version}
\begin{equation}
\label{eq:projection_stabilized}
\psi_h\in M_{h}^{p-1} : \quad \alpha(\psi_h,q)= b(w_h,q)-r(q) \quad \forall q\in M_{h}^{p-1}
\end{equation}
where
\begin{equation*}
\alpha(\psi_h,q) = b(\nabla\psi_h,q) + j_0(\psi_h,q) .
\end{equation*}
Note that this system naturally corresponds to homogeneous Neumann conditions
on $\Gamma_D$ and homogeneous Dirichlet conditions on $\Gamma_N$ (which might
be empty). Now we may define the first projection scheme.

\begin{algorithm}
\label{alg:standard_proj}
The standard projection $\mathcal{P}_h^{\text{std}}$
is given by the following algorithm:
\begin{enumerate}[i)]
\item For any tentative velocity $w_h\in X_h^p$ solve
\begin{equation}
\psi_h\in M_{h}^{p-1}: \quad \alpha(\psi_h,q)= b(w_h,q)-r(q) \quad \forall q\in M_{h}^{p-1}.
\end{equation}
\item Set $\mathcal{P}_h^{\text{std}} w_h = v_h$ where $v_h$ solves
\begin{equation}
\label{eq:std_subtraction}
(v_h,\varphi)_{0,\Omega}  = (w_h,\varphi)_{0,\Omega} - (\nabla\psi_h,\varphi)_{0,\Omega}
\quad \forall \varphi\in X_h^p .
\end{equation}
This requires the solution of a mass matrix which is block-diagonal. Choosing an orthogonal
basis it can even be diagonal and thus the computation is cheap. Note also that this implies
$v_h = w_h-\nabla_h\psi_h$ since $\nabla_h\psi_h\in X_h^p$.
\end{enumerate}
\end{algorithm}
\noindent Unfortunately, this projection is reported to be unstable in the small
time step limit \cite{Steinmoeller2013480} and we also observed this behaviour.
Part of the problem is that $\mathcal{P}_h^{\text{std}}$ is actually not a projection, i.e.
$(\mathcal{P}_h^{\text{std}})^2\neq \mathcal{P}_h^{\text{std}}$.

\subsubsection{Div-div projection}

In order to overcome the stability problem the authors in \cite{2016arXiv160701323K}
suggested to stabilize the projection by an additional term in \eqref{eq:std_subtraction}:

\begin{algorithm}
\label{alg:divdiv_proj}
The div-div projection $\mathcal{P}_h^{\text{div-div}}$
is given by the following algorithm:
\begin{enumerate}[i)]
\item For any tentative velocity $w_h\in X_h^p$ solve (same as before)
\begin{equation*}
\psi_h\in M_{h}^{p-1}: \quad \alpha(\psi_h,q)= b(w_h,q)-r(q) \quad \forall q\in M_{h}^{p-1}.
\end{equation*}
\item Set $\mathcal{P}_h^{\text{div-div}} w_h = v_h$ where $v_h$ solves
\begin{equation}
\label{eq:divdiv_subtraction}
(v_h,\varphi)_{0,\Omega} + \tau_D (\nabla\cdot v_h,\nabla\cdot\varphi)_{0,\Omega}
= (w_h,\varphi)_{0,\Omega} - (\nabla\psi_h,\varphi)_{0,\Omega}
\quad \forall \varphi\in X_h^p
\end{equation}
where $\tau_D$ is a user-supplied constant.
\end{enumerate}
Again this requires the solution of an element-local system
which is not diagonal. As reported in \cite{2016arXiv160701323K} and the examples below
this gives good results with quite small point-wise divergence. However, the projected
velocity does not satisfy a local mass conservation property and
$(\mathcal{P}_h^{\text{div-div}})^2\neq \mathcal{P}_h^{\text{div-div}}$
\end{algorithm}

\subsubsection{Raviart-Thomas projection}

The aim of this subsection is to reconstruct $-\nabla\psi_h$ in the Raviart-Thomas space
of degree $k$ \cite{Brezzi91} on affine cuboid meshes given by
\begin{equation}
\RT_h^k = \{ v\in H(\text{div};\Omega) \mid v|_E\in \RT_E^k \, \forall E\in\mathcal{E}_h\}
\end{equation}
with the Raviart-Thomas space on element $E$ given by
\begin{equation}\label{eq:RTonE}
\RT_E^k = \{ v\in H(\text{div};E) \mid v = T_E(\hat v),
(\hat v)_i = \sum_{\{ \alpha \mid 0\leq\alpha_j\leq k+\delta_{ij}\}} c_{i,\alpha} \hat x^\alpha\}
\end{equation}
where we made use of the Piola transformation of the affine element $E\in\mathcal{E}_h$,
i.e. $\mu_E(\hat x) = B_E \hat x + b_E$, defined as
\begin{equation*}
T_E(\hat v) = \frac{1}{|\det B_E|} B_E \hat v .
\end{equation*}
For $k>0$ the construction needs also the space
\begin{equation}
\Psi_E^k = \{ v\in H(\text{div};E) \mid v = T_E(\hat v),
(\hat v)_i = \sum_{\{ \alpha \mid 0\leq\alpha_j\leq k-\delta_{ij}\}} c_{i,\alpha} \hat x^\alpha\} .
\end{equation}
Note that in contrast to \eqref{eq:RTonE} the polynomial degree in direction $i$ in component $i$
is decreased instead of increased.

Assume that $\psi_h\in M_h^{p-1}$ solves \eqref{eq:projection_stabilized} as before.
Following \cite{ErnHdiv2007} we now compute $\gamma_h=G_h\psi_h \in \RT_h^{k}$, $k=p-1$,
as reconstruction of $-\nabla\psi_h$ as follows. On element $E\in\mathcal{E}_h$ with faces $e\in\partial E$ define
\begin{subequations}
\label{eq:RT_reconstruction}
\begin{align}
(\gamma_h\cdot \normal_e,q)_{0,e} &=
(-\{ \nabla\psi_h \}\cdot \normal_e + \frac{\sigma}{h_e} [\psi_h],q)_{0,e} && e\in\Gamma_h^\text{int}, q\in Q^{k}_e,\\
(\gamma_h\cdot \normal_e,q)_{0,e} &=
(- \nabla\psi_h \cdot \normal_e + \frac{\sigma}{h_e} \psi_h,q)_{0,e} && e\in\Gamma_h^N, q\in Q^{k}_e,\\
(\gamma_h\cdot \normal_e,q)_{0,e} &= 0 && e\in\Gamma_h^D, q\in Q^{k}_e, \label{eq:rt_bc}
\end{align}
and for $k>0$ define in addition
\begin{equation}
\begin{split}
(\gamma_h,r)_{0,E} &= -( \nabla\psi_h,r)_{0,E}
+\frac12 \sum_{e\in\partial E \cap \Gamma_h^\text{int}} (r\cdot \normal_e,[ \psi_h ])_{0,e}
 +\sum_{e\in\partial E \cap \Gamma_h^N} (r\cdot \normal_e, \psi_h)_{0,e} ,
\qquad \forall r \in \Psi_E^k .
\label{eq:RT_projection_c}
\end{split}
\end{equation}
\end{subequations}
With this we can define our final projection method:
\begin{algorithm}
\label{alg:RT_reconstruction}
The RT projection $\mathcal{P}_h^{\text{RT}}$
is given by the following algorithm:
\begin{enumerate}[i)]
\item For any tentative velocity $w_h\in X_h^p$ solve
\begin{equation*}
\psi_h\in M_{h}^{p-1}: \quad \alpha(\psi_h,q)= b(w_h,q)-r(q) \quad \forall q\in M_{h}^{p-1}.
\end{equation*}
\item Reconstruct $\gamma_h=G_h\psi_h \in \RT_h^{p-1}$.
\item Set $\mathcal{P}_h^{\text{RT}} w_h = v_h$ where $v_h$ solves
\begin{equation*}
(v_h,\varphi)_{0,\Omega}  = (w_h,\varphi)_{0,\Omega} + (G_h\psi_h,\varphi)_{0,\Omega}
\quad \forall \varphi\in X_h^p .
\end{equation*}
This requires the solution of a (block-) diagonal system.
\end{enumerate}
\end{algorithm}

\noindent The reconstruction $G_h$ defined above satisfies the following
important property.
\begin{lemma} \label{lemma:projection}
Let $\psi_h\in M_h^{p-1}$ solve $\alpha(\psi_h,q)=l(q)$ for all $q\in M_h^{p-1}$ and
any linear right hand side functional $l$. Let furthermore $\gamma_h=G_h\psi_h \in \RT_h^{p-1}$
be the reconstruction defined above.
Then for every $q\in Q_h^{p-1}$ and $\chi_E$ the characteristic function of
element $E\in\mathcal{E}_h$ we have
\begin{equation}
(\nabla\cdot \gamma_h, q\chi_E)_{0,E} = l(q\chi_E) .
\end{equation}
\begin{proof} Straightforward extension of Theorem 3.1 in \cite{ErnHdiv2007}
from simplicial to affine cuboid elements. Essential ingredients are that for any
$q\in Q_h^{p-1} \Rightarrow \nabla q|_E\in\Psi_E^{p-1}$ and the
special definition of the right hand side in \eqref{eq:RT_projection_c}.
\end{proof}
\end{lemma}
And with this lemma we can prove the following theorem.
\begin{theorem} \label{theo:rt_projection}
The projected velocity $\mathcal{P}_h^{\text{RT}} w_h$ satisfies
the discrete continuity equation exactly, i.e.
\begin{equation}
\label{eq:RT_discrete_continuity}
b(\mathcal{P}_h^{\text{RT}} w_h,q) = r(q) \qquad\forall q\in M_h^{p-1} .
\end{equation}
\begin{proof} The characteristic functions form a partition of unity,
i.e. for any $q\in Q_h^{p-1}$ we have $q=\sum_{E\in\mathcal{E}_h}q\chi_E$. Inserting
into the definition of $b$, observing that
$[\gamma_h]\cdot \normal_e = 0$ since $\gamma_h\in H(\text{div};\Omega)$
as well as $\gamma_h\cdot\normal_e=0$ due to \eqref{eq:rt_bc} and using Lemma \ref{lemma:projection} gives:
\begin{equation}
\begin{split}
b(\mathcal{P}_h^{\text{RT}} w_h,q) &= b(w_h,q) + b(G_h\psi_h,q)\\
&= b(w_h,q) - \sum_{E\in\mathcal{E}_h} (\nabla\cdot \gamma_h,q\chi_E)_{0,E}
+ \sum_{e\in\Gamma_h^\text{int}} ([ \gamma_h ]\cdot \normal_e,\{q\})_{0,e}
 + \sum_{e\in\Gamma_h^D} (\gamma_h\cdot \normal_e, q)_{0,e} \\
&= b(w_h,q) - \sum_{E\in\mathcal{E}_h} l(q\chi_E)
= b(w_h,q) - \sum_{E\in\mathcal{E}_h} [b(w_h,q\chi_E) -r(q\chi_E)]\\
&= b(w_h,q) - b(w_h,q) + r(q)
= r(q)
\end{split}
\end{equation}
\end{proof}
\end{theorem}

\begin{remark} \label{rem:rt_projection} As corollaries we have
\begin{enumerate}[1)]
\item The projected velocity $v_h = \mathcal{P}_h^{\text{RT}} w_h$ satisfies
  the discrete conservation property
  (\ref{eq:local_mass_conservation}) (use the fact \(\chi_E\in
  M_h^{p-1}\) and Theorem \ref{theo:rt_projection}). Note that this
  discrete conservation property can be achieved with reconstruction in
  Raviart-Thomas space with degree \(k\leq p-1\).
\item $(\mathcal{P}_h^{\text{RT}})^2=\mathcal{P}_h^{\text{RT}}$ follows
from Theorem \ref{theo:rt_projection} and the fact
  that $l(q) = b(\mathcal{P}_h^{\text{RT}} w_h,q) - r(q) =0, q\in M_h^{p-1},$
  is the right-hand side in step i) of Algorithm
  \ref{alg:RT_reconstruction}. Therefore when applying
  \(\mathcal{P}_h^{\text{RT}}\) twice a zero correction is produced in
  the second application.
\end{enumerate}
\end{remark}

The discrete continuity equation does not imply that the divergence of the
projected velocity vanishes point-wise. The following Lemma shows that
the divergence in the interior of elements is controlled in an integral sense only by the
jumps of the tentative velocity:
\begin{lemma} \label{lemma:rt_pointwise_divergence}
The projected velocity $v_h=\mathcal{P}_h^{\text{RT}} w_h$ satisfies
for all $q\in M_h^{p-1}$, $E\in\mathcal{E}_h$ and $q_E=q\chi_E$:
\begin{equation}
(\nabla\cdot v_h,q_E)_{0,E} =
\frac12 \sum_{e\in\Gamma_h^\text{int}\cap\partial E} ([ w_h ]\cdot\normal_e,q_E)_{0,e}
+ \sum_{e\in\Gamma_h^{D}\cap\partial E} ( (w_h-g)\cdot \normal_e,q_E )_{0,e} .
\end{equation}
\begin{proof} Using Lemma \ref{lemma:projection} we get
\begin{equation}
\begin{split}
(\nabla \cdot v_h,q_E)_{0,E} &= (\nabla \cdot w_h,q_E)_{0,E} + (\nabla \cdot \gamma_h,q_E)_{0,E}
= (\nabla \cdot w_h,q_E)_{0,E} + l(q_E) \\
&= (\nabla \cdot w_h,q_E)_{0,E} + b_h(w_h,q_E)-r(q_E)\\
&= (\nabla \cdot w_h,q_E)_{0,E} - (\nabla\cdot w_h, q_E)_{0,E}
+\frac12 \sum_{e\in\Gamma_h^\text{int}\cap\partial E} ([ w_h ]\cdot \normal_e,q_E)_{0,e}
+ \sum_{e\in\Gamma_h^{D}\cap\partial E} ( (w_h-g)\cdot \normal_e,q_E )_{0,e} .
\end{split}
\end{equation}
\end{proof}
\end{lemma}

\subsection{Pressure-correction schemes}
\label{sec:press-corr-schem}

Since the nonlinear term in the Navier-Stokes equations does not play an
essential role in the derivation of the projection methods we hereafter
consider the instationary Stokes equations. The equations in the
subproblems arise from the method of lines discretization.

\subsubsection{Incremental pressure-correction scheme (IPCS)}
\label{sec:incr-press-corr}

The IPCS is a straightforward way to split between incompressibility and
dynamics. In the viscous substep the pressure is made explicit that we
denote by \(p_h^{\bigstar, k+1}\). In the second substep a pressure
correction is computed to accordingly correct the velocity. The
particular choice of the time discretization is not important. It is
possible to use the implicit Euler time stepping or second order time
stepping methods such as BDF2 or Alexander's second order strongly
S-stable scheme \cite{Alexander:1977:DIR}. The semi-discretized in space
splitting scheme then reads as follows:

\begin{enumerate}
\item Tentative velocity step, compute $\tilde{v}_h^{k+1}$:
  \begin{equation*}
    \rho(\partial_t v_h, \varphi_h) + \mu a(v_h,\varphi_h) +
    b(\varphi_h,p_h^{\bigstar, k+1}) = l(\varphi_h;t) \quad \forall
    \varphi_h\in X_h^p
  \end{equation*}
\item Projection step: Compute \(\delta p_h^{k+1}
    = 1 / \Delta t^{k+1} \psi_h\) and \(v_h^{k+1} =
    \mathcal{P}_h\tilde{v}_h^{k+1}\) by choosing one of the projectors
    given by Algorithm \ref{alg:standard_proj}, \ref{alg:divdiv_proj} or
    \ref{alg:RT_reconstruction}.
\item Pressure update:
  \begin{equation*}
    p_h^{k+1} = p_h^{\bigstar, k+1} + \delta p_h^{k+1} .
  \end{equation*}
\end{enumerate}
The choice \(p_h^{\bigstar, k+1} = 0\), implicit Euler as time stepping
yields to Chorin's projection method. Constant extrapolation
\(p_h^{\bigstar, k+1} = p_h^k\) gives the IPCS.
The IPCS introduces the artificial boundary conditions for the
pressure correction which lead to the series of equalities
\begin{align}
  \left . \partial_{\normal} p_h^{k+1} \right|_{\Gamma_D} &= \hdots =
  \left . \partial_{\normal} p_h^1 \right|_{\Gamma_D} = \left
    . \partial_{\normal} p_h^0 \right|_{\Gamma_D} \\
  \left . p_h^{k+1} \right|_{\Gamma_N} &= \hdots = \left . p_h^1
  \right|_{\Gamma_N} = \left . p_h^0 \right|_{\Gamma_N}
\end{align}
for the pressure itself over time. In the purely Dirichlet case,
i.e. $\Gamma_N = \emptyset$, the scheme is fully first-order accurate
even if the implicit Euler time stepping is used. But when $\Gamma_N
\neq \emptyset$ the order of approximation of the velocity in the
$H_0^1$-norm and of the pressure in the $L^2$-norm is degraded due to
the homogeneous Dirichlet boundary conditions for the pressure.

There is little improvement regarding the order of the scheme when a
second order time stepping method is used. In the purely Dirichlet case
the scheme is fully second order on the velocity in the $L^2$-norm but
it stays first order on the velocity in the $H_0^1$-norm and on the
pressure in the $L^2$-norm. For $\Gamma_N \neq \emptyset$ the
approximation order even stays the same.

The constant extrapolation for the explicit pressure in the momentum
equation implies that the scheme has an irreducible splitting error of
$\mathcal{O}(\Delta t^2)$. Hence using a higher than second order time
discretization does not improve the overall accuracy.

\subsubsection{Rotational incremental pressure-correction scheme (RIPCS)}
\label{sec:rotat-incr-press}

One reason for the above scheme to have poor convergence properties
especially when outflow boundary conditions are present is that the
pressure boundary conditions stay constant over time. To overcome this
difficulty it was first introduced by Timmermans, Minev and Van De Vosse
\cite{FLD:FLD373} to use the \emph{rotational form} of the Laplacian, namely
\begin{equation}
  \label{eq:laplacian_rotational_form}
  -\Delta v = \nabla\times(\nabla\times v) - \nabla(\nabla\cdot v) .
\end{equation}
To understand why this modification performs better we consider for
simplicity the momentum equation in classical form and insert the
rotational form of the Laplacian:
\begin{equation}
  \frac{\tilde{v}_h^{k+1} - v_h^k}{\Delta t^{k+1}} +
  \mu\nabla\times(\nabla\times\tilde{v}_h^{k+1}) +
  \nabla(p_h^{\bigstar, k+1} - \mu\nabla\cdot\tilde{v}_h^{k+1}) =
  f(t^{k+1})
\end{equation}
where $p_h^{\bigstar, k+1}$ is as before an approximation of
$p(t^{k+1})$. Eliminating the tentative velocity
\(\tilde{v}_h^{k+1} = v_h^{k+1} + \Delta t^{k+1}\nabla\delta p_h^{k+1}\)
with the Helmholtz decomposition gives
\begin{equation}
  \frac{v_h^{k+1} - v_h^k}{\Delta t^{k+1}} +
  \mu\nabla\times(\nabla\times\tilde{v}_h^{k+1}) + \nabla(\delta p_h^{k+1}
  + p_h^{\bigstar, k+1} - \mu\nabla\cdot\tilde{v}_h^{k+1}) = f(t^{k+1}) \; .
\end{equation}
Thus the quantity $\delta p_h^{k+1} + p_h^{\bigstar, k+1} -
\mu\nabla\cdot\tilde{v}_h^{k+1}$ can be interpreted as an approximation
of the pressure. Hence retaining the time step with the momentum
equation the tables can be turned to obtain the \emph{incremental
  pressure-correction scheme in rotational form}:

\begin{enumerate}
\item Tentative velocity step, compute $\tilde{v}_h^{k+1}$:
  \begin{equation*}
    \rho(\partial_t v_h, \varphi_h) + \mu a(v_h,\varphi_h) +
    b(\varphi_h,p_h^{\bigstar, k+1}) = l(\varphi_h;t) \quad \forall
    \varphi_h\in X_h^p
  \end{equation*}
\item Projection step: Compute \(\delta p_h^{k+1}
    = 1 / \Delta t^{k+1} \psi_h\) and \(v_h^{k+1} =
    \mathcal{P}_h\tilde{v}_h^{k+1}\) by choosing one of the projectors
    given by Algorithm \ref{alg:standard_proj}, \ref{alg:divdiv_proj} or
    \ref{alg:RT_reconstruction}.
\item Pressure update with scaling factor $\omega$:
  \begin{equation*}
    (p_h^{k+1}, q_h) = (\omega\delta p_h^{k+1} + p_h^{\bigstar, k+1}, q_h) + \mu
    (b(\tilde{v}_h^{k+1}, q_h) - r(q_h;t^{k+1})) \quad \forall q_h\in M_h^{p-1} \; .
  \end{equation*}
\end{enumerate}
The scaling factor is usually set to $\omega = 1$ for first order time
stepping schemes and to $\omega = \frac 32$ for second order time
stepping schemes.

The contribution $\nabla\cdot\tilde{v}_h^{k+1}$ improves the accuracy of
the scheme such that it is first order accurate for both Dirichlet and
outflow boundary conditions. The use of a second order time stepping
scheme improves the convergence rate on the velocity in the $H_0^1$-norm
and on the pressure in the $L^2$-norm to $\frac32$ when $\Gamma_N =
\emptyset$. In the presence of outflow boundary conditions the
convergence rate $\frac32$ for the velocity in the $L^2$-norm is likely
to be the best possible whereas the convergence rate in the $H_0^1$-norm
for the velocity and in the $L^2$-norm for the pressure is limited to 1.
As in the IPCS higher than second order time stepping schemes do not
improve the overall accuracy.

\section{Numerical experiments}
\label{sec:numer-exper}

We start the numerical experiments by cross-comparing the pointwise
divergence and local mass conservation for the div-div projection and
the $H(\text{div})$ reconstruction. Then we illustrate the convergence properties
of the IPCS and RIPCS for
global Dirichlet boundary conditions \ref{sec:glob-dirichl-bound}, mixed
boundary conditions \ref{sec:mixed-bound-cond}, periodic boundary
conditions \ref{sec:peri-bound-cond} and also in 3D using the Beltrami
flow problem \ref{sec:beltrami-flow}. Both schemes are tested in their
second order formulation.
Temporal convergence is analyzed for the Taylor-Hood-like DG-spaces
\THcube{2}{1}, \THcube{3}{2}, \THcube{4}{3} and also local mass
conservation - given as the left-hand side of
(\ref{eq:local_mass_conservation}) - is investigated.

\subsection{Solvers and Implementation}
\label{sec:solv-impl}

The parallel solver has been implemented in a high-performance C++ code
based on the DUNE discretization framework \cite{dune08:1,Bastian2016}.
The assembly of residuals and jacobians uses spectral discontinuous
Galerkin methods. Sum-factorization technique for tensor product bases
is employed that reduce the computational complexity significantly.
Every velocity component underlies the same ansatz
space. Therefore sum-factorization applied to the scalar
convection-diffusion equation as described in \cite{Kronbichler2012} can be
expanded in a straightforward way to the subproblems in the splitting
schemes. The viscous substep is solved with a matrix-free Newton method
with a single block SOR preconditioner in GMRes as a linear solver. The
pressure Poisson equation is solved with hybrid AMG-DG preconditioner
where the correction in the conforming \(\mathcal{Q}_1\) subspace is
rediscretized, \cite{amg4dg} and the matrix
on the DG level is not required for this purpose. Thus it is possible to
do either matrix-free or matrix-based operator application and smoothing
on the DG level.

\subsection{Local mass conservation}
\label{sec:local-mass-cons-1}

We consider the Navier-Stokes equations on the domain $\Omega =
(-1,1)^2$ and take the two dimensional
Taylor-Green vortex which has been studied before by
\cite{10.2307/96892,ChorinProjectionNavierStokes,5982828}.
In two dimensions the Taylor-Green vortex possesses the exact solution
\begin{align}
  v_1(x,y,t) &= -e^ {- 2\,\pi^2\,\nu\,t }\,\cos
  \left(\pi\,x\right)\,\sin \left(\pi\,y\right) \notag \\
  v_2(x,y,t) &= e^ {- 2\,\pi^2\,\nu\,t }\,\sin \left(\pi\,x\right)\,
  \cos \left(\pi\,y\right) \notag \\
  p(x,y,t) &= -0.25\,\rho\,e^ {- 4\,\pi^2\,\nu\,t }\,\left(\cos
    \left(2\,\pi\,y \right)+\cos
    \left(2\,\pi\,x\right)\right) \label{eq:periodic_testproblem} \; .
\end{align}
The source term is given by $f = 0$. We set $\rho = 1$, $\mu = 1/100$
and $\nu = \mu / \rho$. Periodic boundary conditions are
imposed in both the $x$ and $y$ directions.
We do the computations on a $160 \times 160$ rectangular mesh. The discussion
on the temporal convergence rates is postponed to Section
\ref{sec:peri-bound-cond}.

\begin{figure}[!htpb]
  \includegraphics[width=0.5\textwidth]{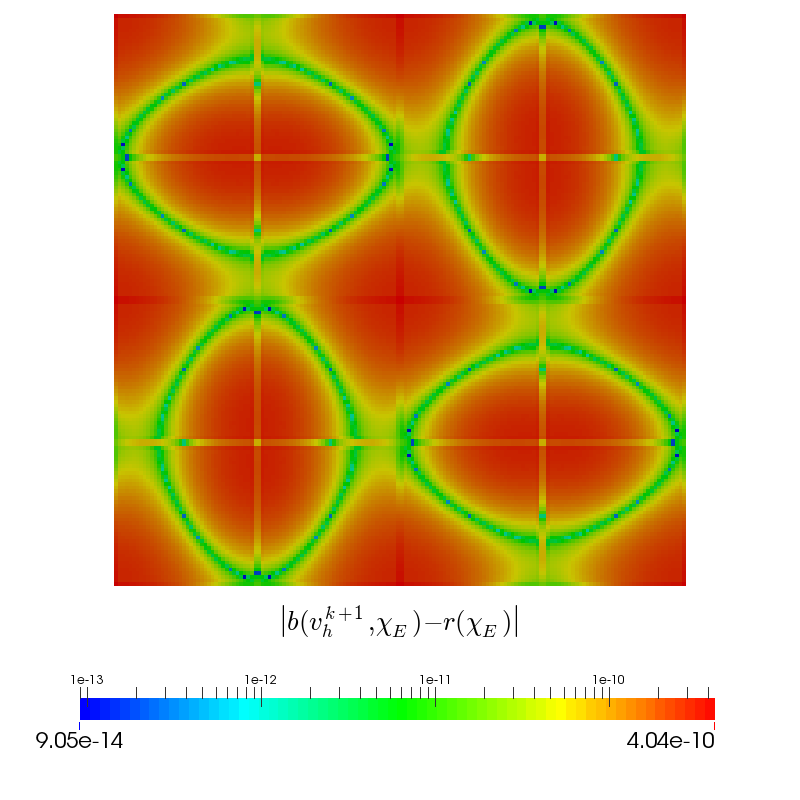}
  \includegraphics[width=0.5\textwidth]{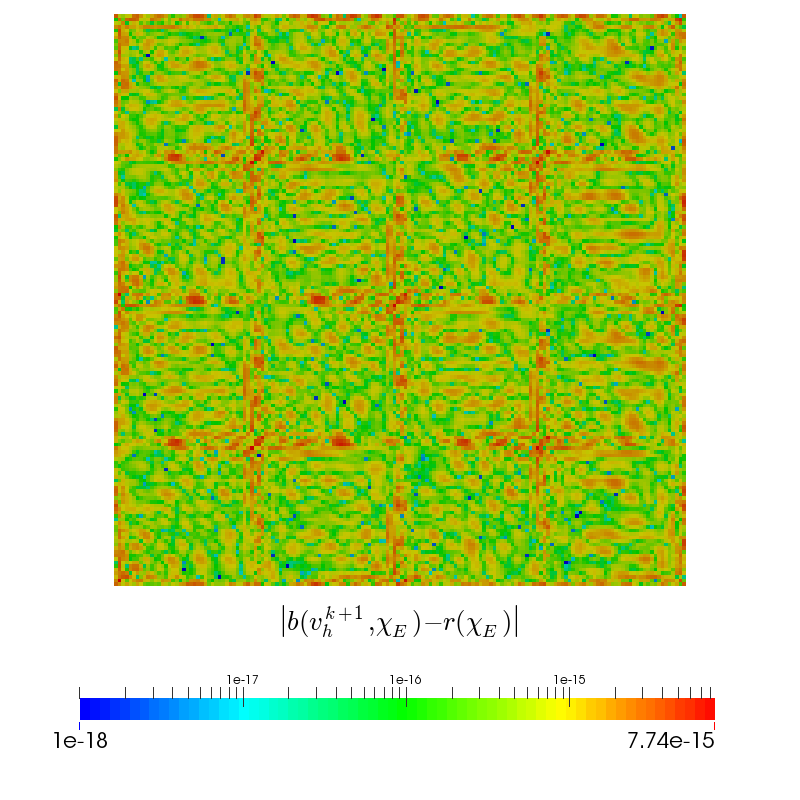}
  \caption{Local mass conservation of the Taylor-Green vortex solution
    at time 1 with \(\Delta t = 0.025\) obtained by the RIPCS. \newline
    Left part shows \THcube{2}{1} with div-div projection.
    Right part shows \(\mathcal{Q}_2 / \mathcal{Q}_1\) with
    reconstruction \(\RT_h^0\), identical to reconstruction in \(\RT_h^1\).}
  \label{pic:ripcs_order2_mass_conservation}
\end{figure}

We start the discussion on the choice of order in the Raviart-Thomas
space. We have shown in Theorem \ref{theo:rt_projection} that for
\(\RT_h^{p-1}\) it holds:
(I) \((\mathcal{P}_h^{\text{RT}})^2=\mathcal{P}_h^{\text{RT}}\),
(II) the reconstructed velocity satisfies the continuity equation and
(III) is locally mass conservative.
However a naive approach by looking at the dimension of the local
function space of \(\nabla_h M_h^{p-1}\) also accounts to possibly choose
\(\RT_h^{p-2}\). As stated in Remark \ref{rem:rt_projection} local
mass conservation can still be achieved with reconstruction in
Raviart-Thomas space of degree \(p-2\). This is demonstrated on the
right of figure \ref{pic:ripcs_order2_mass_conservation} and notably we
get the same distribution with \(\RT_h^{p-1}\). Moreover numerical
experiments with the power iteration applied to the operator
\(\mathcal{P}_h^{\text{RT}}\) have shown that
\((\mathcal{P}_h^{\text{RT}})^2 \tilde{v}_h^{k+1} =
\mathcal{P}_h^{\text{RT}} \tilde{v}_h^{k+1}\) also for
\(\RT_h^{p-2}\). Table \ref{tab:ripcs_hdiv_order2_periodic_rates}
- \ref{tab:ripcs_rt1_order2_periodic_rates} compare the temporal
accuracy between the discretizations
  \THcube{2}{1} with reconstruction in \(\RT_h^0\) or
  \(\RT_h^1\). It can be seen that there is no significant
difference on the error at final time. Reconstruction in the
\(\RT_h^{p-2}\) space provides thus to be a sufficient alternative
in the splitting algorithm.

Next we want to cross-compare the temporal accuracy for the spatial
discretizations \THcube{2}{1}, \THcube{3}{2} with the div-div projection
and \THcube{2}{1} with reconstruction in \(\RT_h^0\), \THcube{3}{2} with
reconstruction in \(\RT_h^1\).
Table \ref{tab:ripcs_order2_periodic_rates} -
\ref{tab:ripcs_hdiv_order2_periodic_rates} show the errors for the RIPCS
\THcube{2}{1} with div-div projection and the RIPCS
\THcube{2}{1} with reconstruction in \(\RT_h^0\)
and table \ref{tab:ripcs_order3_periodic_rates} -
\ref{tab:ripcs_hdiv_order3_periodic_rates} the errors for the RIPCS
\THcube{3}{2} with div-div projection and the RIPCS
\THcube{3}{2} with reconstruction in \(\RT_h^1\), respectively.
There is no significant difference in the temporal
behaviour for both pairs, a logarithmic plot of the errors would
lead to indistinguishable curves.
Thus for the upcoming investigation on
the convergence properties we will use the div-div projection
technique because it is an inexpensive alternative to the $H(\text{div})$
reconstruction which is at the time only implemented up to order one.
Note that the errors in the tables
\ref{tab:ripcs_order2_periodic_rates},
\ref{tab:ripcs_order3_periodic_rates} are also contained in the figures
of \ref{fig:ipcs_ripcs_periodic_rates}.

\begin{table}[H]
  \centering
  \input{convergence-tables/ripcs_order2_periodic_rates.out.tex}
  \caption{Errors for the Taylor-Green vortex at
    final time T=2 obtained by RIPCS and \THcube{2}{1} with div-div
    projection}
  \label{tab:ripcs_order2_periodic_rates}
\end{table}

\begin{table}[H]
  \centering
  \input{convergence-tables/ripcs_hdiv_order2_periodic_rates.out.tex}
  \caption{Errors for the Taylor-Green vortex at
    final time T=2 obtained by RIPCS and \(\mathcal{Q}_2 / \mathcal{Q}_1\)
    with reconstruction in \(\RT_h^0\)}
  \label{tab:ripcs_hdiv_order2_periodic_rates}
\end{table}

\begin{table}[H]
  \centering
  \input{convergence-tables/ripcs_rt1_order2_periodic_rates.out.tex}
  \caption{Errors for the Taylor-Green vortex at final time T=2 obtained
    by RIPCS and \THcube{2}{1} with
      reconstruction in \(\RT_h^1\)}
  \label{tab:ripcs_rt1_order2_periodic_rates}
\end{table}

\begin{table}[H]
  \centering
  \input{convergence-tables/ripcs_order3_periodic_rates.out.tex}
  \caption{Errors for the Taylor-Green vortex at
    final time T=2 obtained by RIPCS and \THcube{3}{2} with div-div
    projection}
  \label{tab:ripcs_order3_periodic_rates}
\end{table}

\begin{table}[H]
  \centering
  \input{convergence-tables/ripcs_hdiv_order3_periodic_rates.out.tex}
  \caption{Errors for the Taylor-Green vortex at
    final time T=2 obtained by RIPCS and
    \THcube{3}{2} with reconstruction in \(\RT_h^1\)}
  \label{tab:ripcs_hdiv_order3_periodic_rates}
\end{table}

\begin{figure}[!htpb]
  \includegraphics[width=0.5\textwidth]{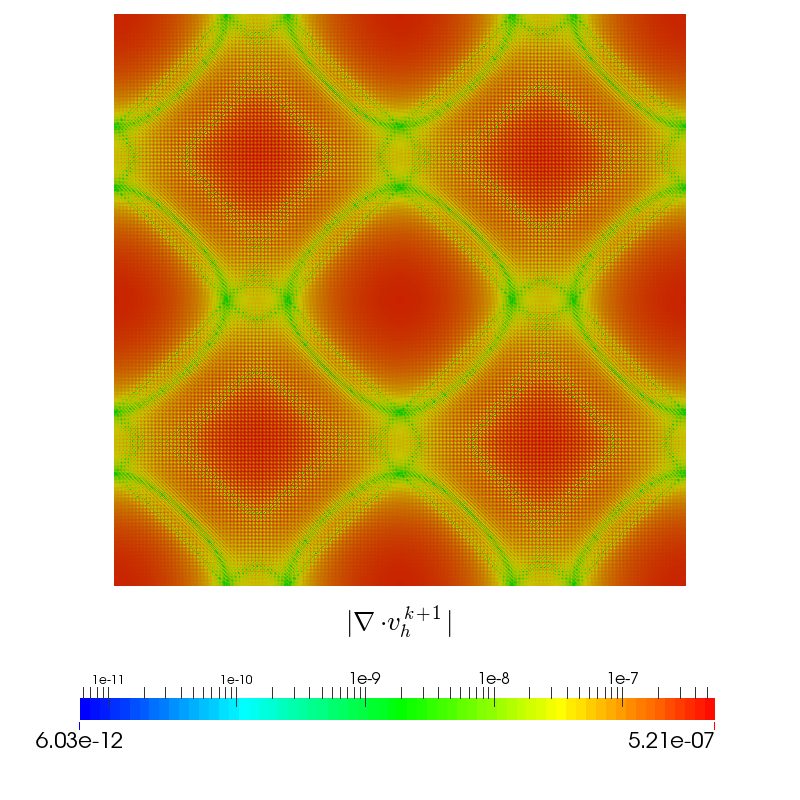}
  \includegraphics[width=0.5\textwidth]{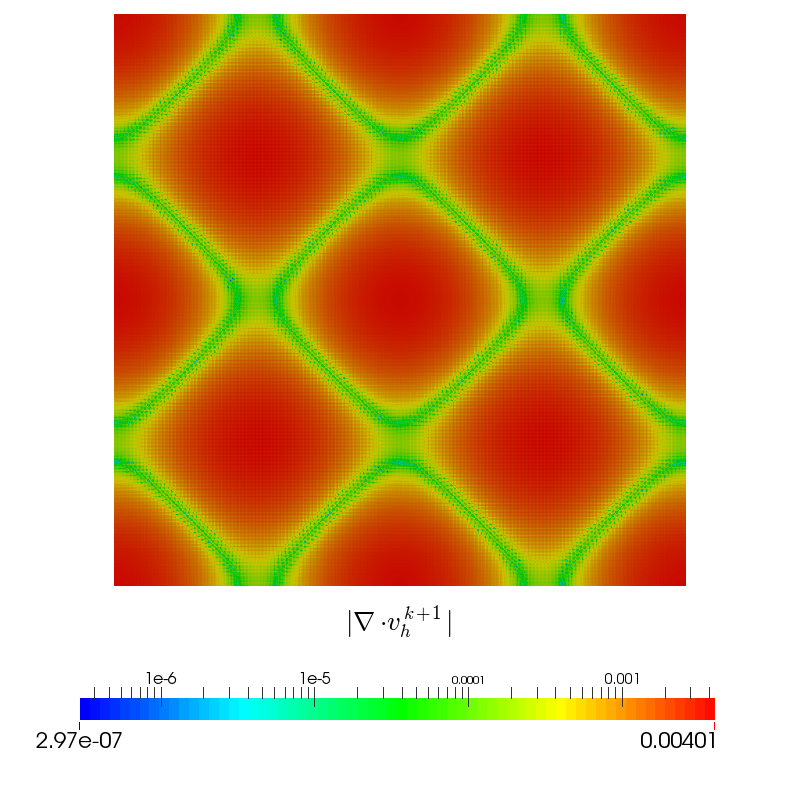}
  \caption{Pointwise divergence of the Taylor-Green vortex solution at
    time 1 with \(\Delta t = 0.025\) obtained by the RIPCS. \newline
    Left part shows \THcube{2}{1} with div-div projection.
    Right part shows \(\mathcal{Q}_2 / \mathcal{Q}_1\) with
    reconstruction in \(\RT_h^0\).}
  \label{pic:ripcs_order2_pointwise_div}
\end{figure}

In figure \ref{pic:ripcs_order2_pointwise_div} the pointwise divergence
for \(p=2\) on each mesh element is presented. The
element-local div-div projection leads to smaller pointwise
divergence than obtained with the $H(\text{div})$
reconstruction. But it does not really cure the error on the local mass
conservation. Compared to the standard \(L^2\)-projection the
div-div projection reduces the values of the pointwise divergence
and local mass conservation. The magnitude of the pointwise divergence
from the $H(\text{div})$ reconstruction is in between the magnitudes
from the standard \(L^2\)-projection and the
  stabilized variant, it is not identically zero as
predicted by Lemma \ref{lemma:rt_pointwise_divergence}. The distribution
of the divergence error with \THcube{2}{1} and
reconstruction in \(\RT_h^1\) is similar and has the
same maximum.

Figure \ref{pic:ripcs_order2_mass_conservation} shows the error on local
mass conservation for \(p=2\). According to our discussion at the
beginning of \ref{sec:local-mass-cons-1} this appealing conservation
property is perfectly fulfilled for the \(\RT_h^{p-1}\) and
\(\RT_h^{p-2}\) reconstructions of the Helmholtz correction.

\subsection{Global Dirichlet boundary conditions}
\label{sec:glob-dirichl-bound}

We consider the Stokes equations on the domain $\Omega = (0,1)^2$ and
take the exact solution to be
\begin{align}
  v_1(x,y,t) &= \sin \left(x+t\right) \,\sin \left(y+t\right) \notag \\
  v_2(x,y,t) &= \cos \left(x+t\right) \,\cos \left(y+t\right) \notag \\
  p(x,y,t) &= \sin
  \left(x-y+t\right) \label{eq:nonhom_dirichlet_testproblem} \; .
\end{align}
The source term is given by $f = \partial_t v - \Delta v + \nabla
p$. The density and viscosity are both set to $\rho = \mu =
1$. Computations were done on a \(160\times 160\) rectangular mesh.

\begin{figure}[!htpb]
  \includegraphics[width=0.5\textwidth]{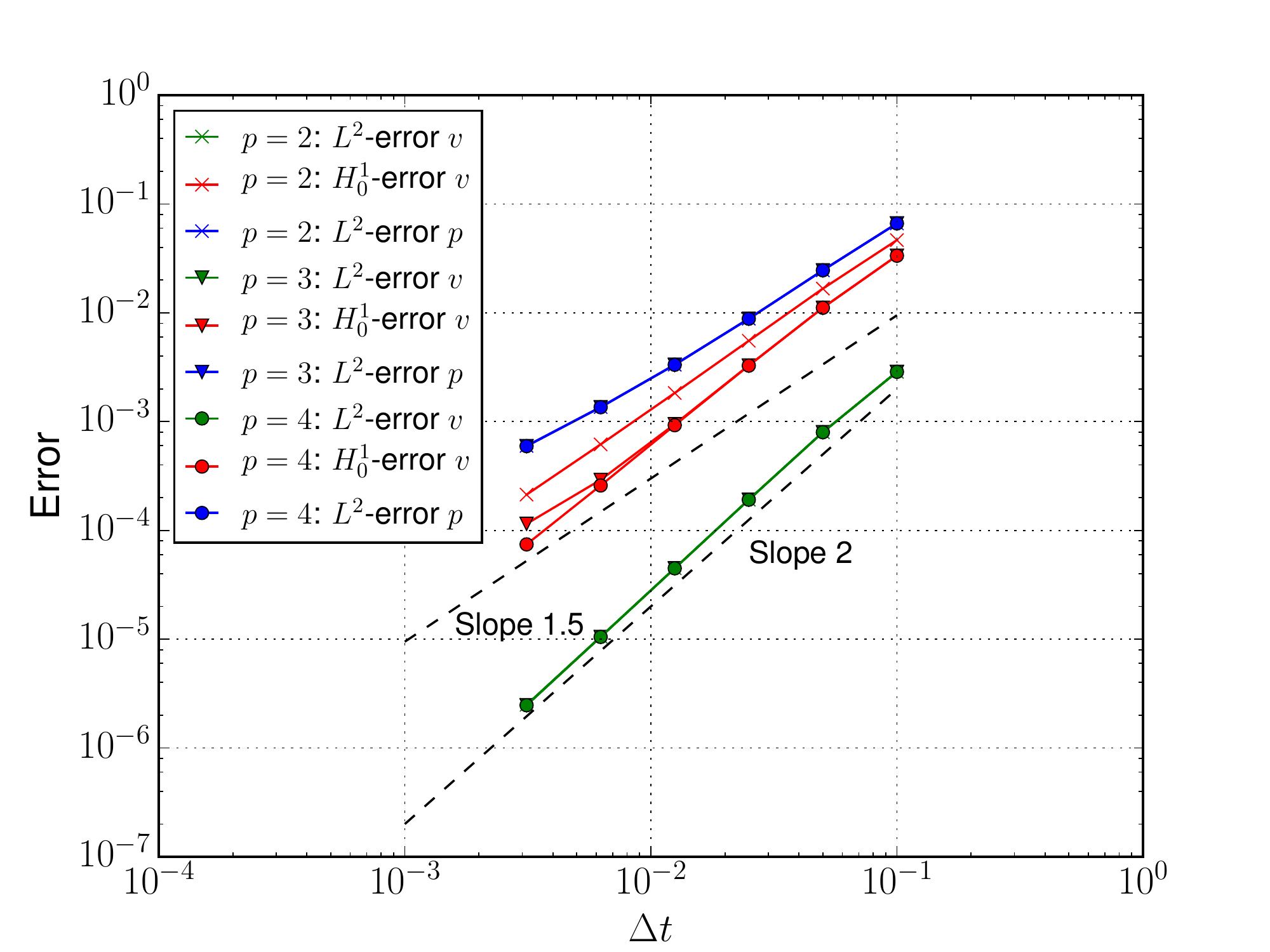}
  \includegraphics[width=0.5\textwidth]{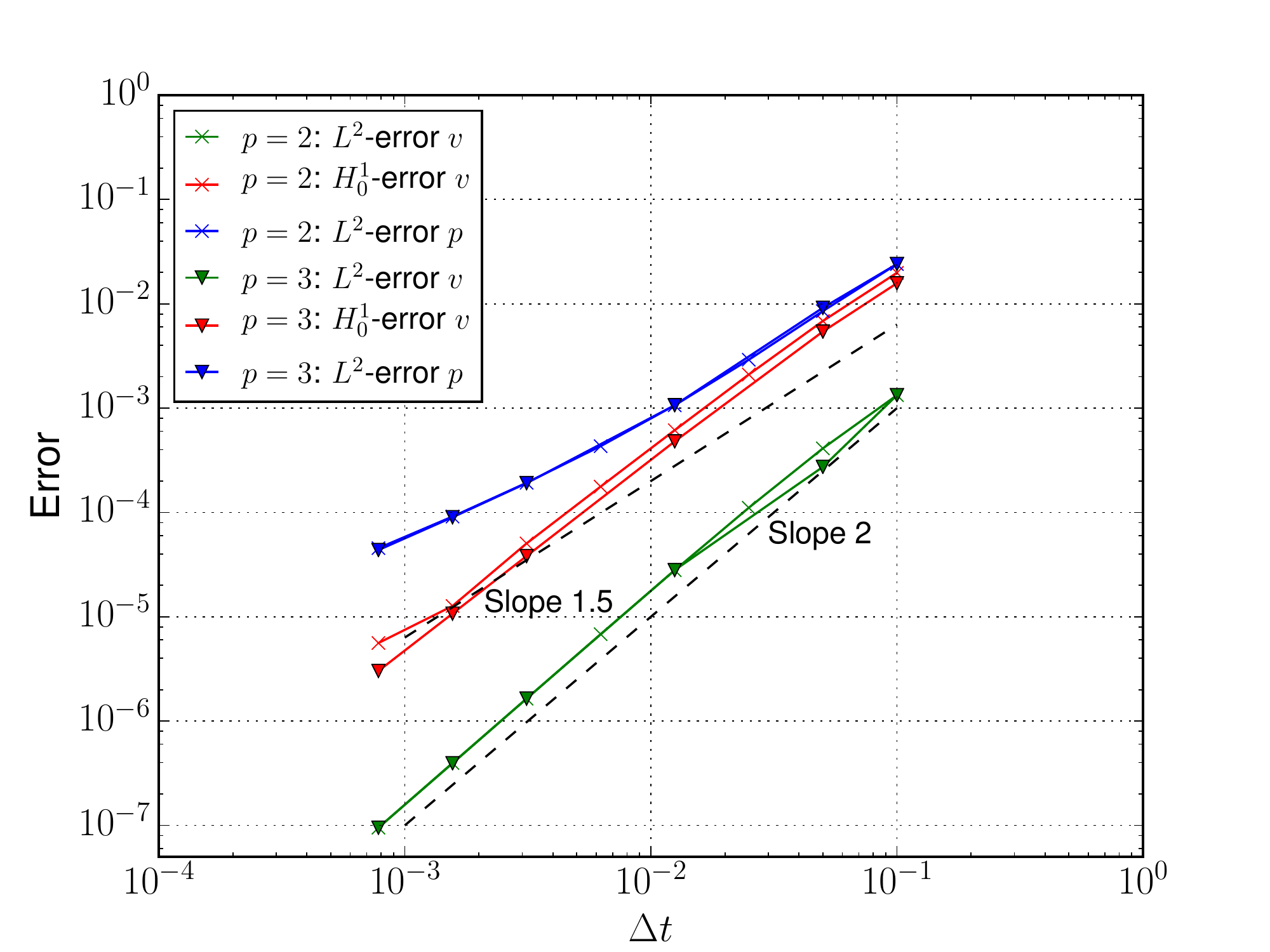}
  \caption{Errors and convergence rates at final time T=1 for the global
    Dirichlet problem and spatial discretizations \THcube{2}{1},
    \THcube{3}{2}, \THcube{4}{3}.
    Left part shows the IPCS.
    Right part shows the RIPCS.}
  \label{fig:ipcs_ripcs_dirichlet_rates}
\end{figure}

Figure \ref{fig:ipcs_ripcs_dirichlet_rates} shows the error and the
convergence rates as function of \(\Delta t\) for the IPCS and
RIPCS. The green curves show the \(L^2\)-error for the velocity, the red
curves the \(H_0^1\)-error for the velocity and the blue curves the
\(L^2\)-error for the pressure obtained by the polynomial degrees
\(p=2,3,4\). The curves grouped by the same color are almost identical
meaning that the splitting error is dominant in the measured range of
\(\Delta t\). Therefore we have left out the curves with \(p=4\) on the
right for the sake of clarity. A transition towards smaller time steps
causes earlier flattening out of the error curves the lower the spatial
order is. This emerges at first for the \(H_0^1\)-error for the velocity
and \(L^2\)-error for the pressure. This is demonstrated for the
Taylor-Green vortex solution in section \ref{sec:peri-bound-cond},
c.f. right of figure \ref{fig:ipcs_ripcs_periodic_rates}.

Theory states that the solution of the second order IPCS satisfies the
following error estimates:
(I) $L^2$-velocity: $\mathcal{O}(\Delta t^2)$
(II) $H_0^1$-velocity, $L^2$-pressure: $\mathcal{O}(\Delta t)$.
On the left of figure \ref{fig:ipcs_ripcs_dirichlet_rates} it is
observed that the velocity error in the \(L^2\)-norm is second order
accurate, in the other two error measures the rate is 1.5 which is
better than the prediction. Now the solution of the RIPCS satisfies the
following error estimates:
(I) $L^2$-velocity: $\mathcal{O}(\Delta t^2)$
(II) $H_0^1$-velocity, $L^2$-pressure: $\mathcal{O}(\Delta t^{\frac32})$.
The convergence rates on the right of figure
\ref{fig:ipcs_ripcs_dirichlet_rates} are consistent with the error
estimates. Note that the \(L^2\)-errors on the velocity and pressure are
almost identical to the results presented in \cite{Guermond04onthe}. The
reason for the slight difference is likely to be the usage of BDF2 in
\cite{Guermond04onthe} as time stepping.

A consideration of local mass conservation shows that it is well satisfied
in the interior of the domain for the div-div projection. However the
largest values \(\sim 10^{-9}\) are located in the cells that share an
edge with the boundary. This is due to the artificial boundary
conditions on the pressure.

The situation is different for the $H(\text{div})$ reconstruction. In that case
the distribution is similar to the right in figure
\ref{pic:ripcs_order2_mass_conservation} with \(\max_E
|b(v_h^{k+1},\chi_E) - r(\chi_E)| \sim 5\cdot 10^{-14}\).

\subsection{Mixed boundary conditions}
\label{sec:mixed-bound-cond}

We consider again the Stokes equations on the domain $\Omega = (0,1)^2$
and take the exact solution to be
\begin{align}
  v_1(x,y,t) &= \sin x\,\sin \left(y+t\right) \notag \\
  v_2(x,y,t) &= \cos x\,\cos \left(y+t\right) \notag \\
  p(x,y,t) &= \cos x\,\sin
  \left(y+t\right) \label{eq:outflow_testproblem} \; .
\end{align}
The source term is again given by $f = \partial_t v - \Delta v + \nabla
p$. The density and viscosity are both set to \(\rho = \mu =
1\). The outflow boundary is located at $\Gamma_N =
\{(x,y)\in\partial\Omega \mid x=0 \}$.
Computations were done on \(160\times 160\) rectangular mesh.

\begin{figure}[!htpb]
  \includegraphics[width=0.5\textwidth]{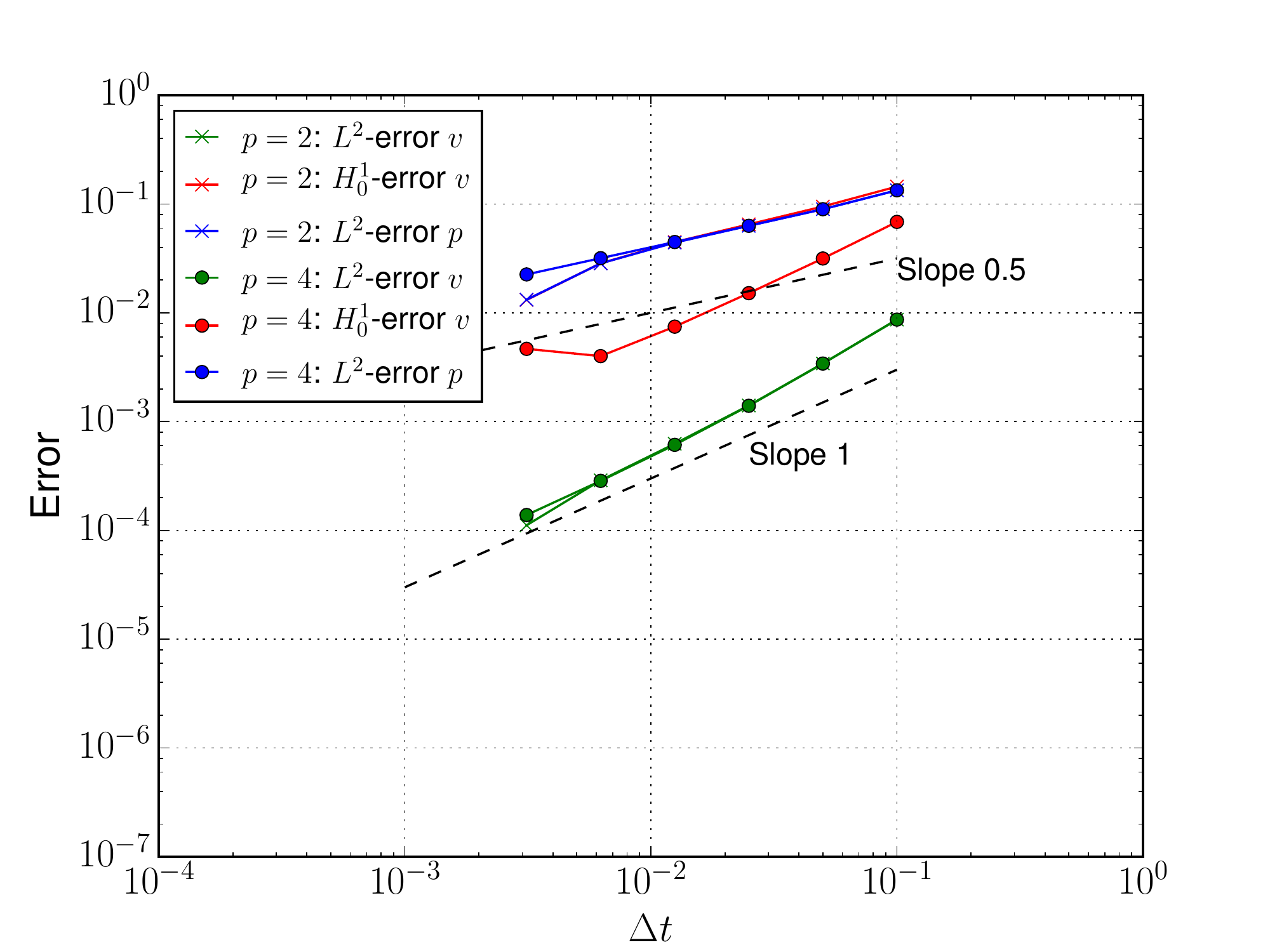}
  \includegraphics[width=0.5\textwidth]{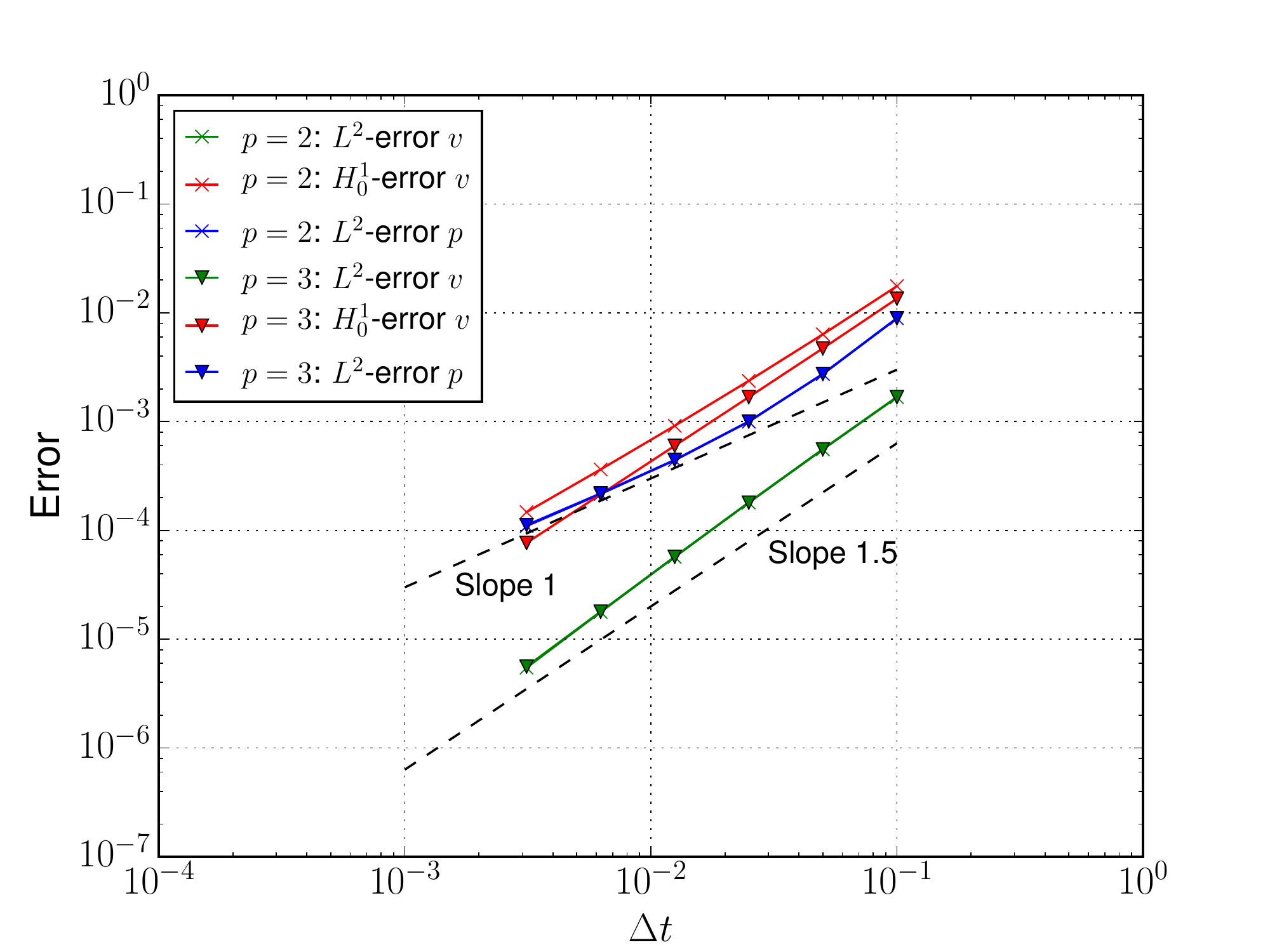}
  \caption{Errors and convergence rates at final time T=1 for the mixed
    boundary condition problem and spatial discretizations
    \THcube{2}{1}, \THcube{3}{2}, \THcube{4}{3}.
    Left part shows the IPCS.
    Right part shows the RIPCS.}
  \label{fig:ipcs_ripcs_outflow_rates}
\end{figure}

Figure \ref{fig:ipcs_ripcs_outflow_rates} shows the error and the
convergence rates as function of \(\Delta t\) for the IPCS and
RIPCS. The green curves show the \(L^2\)-error for the velocity, the red
curves the \(H_0^1\)-error for the velocity and the blue curves the
\(L^2\)-error for the pressure obtained by the polynomial degrees
\(p=2,3,4\). The curves grouped by the same color are almost identical
meaning that the splitting error is dominant in the measured range of
\(\Delta t\). Therefore we have left out the curves with \(p=3\) on the
left and \(p=4\) on the right for the sake of clarity. A transition
towards smaller time steps causes earlier flattening out of the error
curves the lower the spatial order is. This emerges at first for the
\(H_0^1\)-error for the velocity and \(L^2\)-error for the
pressure. This is demonstrated for the Taylor-Green vortex solution in
section \ref{sec:peri-bound-cond}, c.f. right of figure
\ref{fig:ipcs_ripcs_periodic_rates}.

The solution of the IPCS satisfies the following error estimates:
(I) $L^2$-velocity: $\mathcal{O}(\Delta t)$
(II) $H_0^1$-velocity, $L^2$-pressure: $\mathcal{O}(\Delta t^{\frac12})$
which are identical to the first order IPCS. The results on the left of
figure \ref{fig:ipcs_ripcs_outflow_rates} indeed show that the pressure
approximation is poor due to the homogeneous Dirichlet boundary
condition imposed on $\Gamma_N$.
The RIPCS delivers improved error estimates in presence of mixed
boundary conditions:
(I) $L^2$-velocity: $\mathcal{O}(\Delta t^{\frac32})$
(II) $H_0^1$-velocity, $L^2$-pressure: $\mathcal{O}(\Delta t)$.
The convergence rates on right of figure
\ref{fig:ipcs_ripcs_outflow_rates} are consistent with those
estimates. Furthermore the error on the velocity in the $L^2$-norm
behaves like $\mathcal{O}(\Delta t^{\frac53})$ which
is also observed in Guermond, Minev and Shen
\cite{Guermond06,doi:10.1137/040604418}.
The error in the $H_0^1$-norm is close to $\mathcal{O}(\Delta
t^{\frac54})$ which is higher than the rate $\mathcal{O}(\Delta t)$
predicted by theory.
Note that \cite{Guermond06,doi:10.1137/040604418} have used BDF2 as time
stepping for this problem and therefore the error curves are almost
identical.

Another consideration of local mass conservation shows that it is well
satisfied in the interior of the domain. But due to the homogeneous
Dirichlet boundary conditions for the pressure imposed on \(\Gamma_N\),
the largest errors \(\sim 10^{-7}\) are located in the cells next to
outflow boundary.

With the $H(\text{div})$ reconstruction we have
\(\max_E |b(v_h^{k+1},\chi_E) - r(\chi_E)|
\sim 5\cdot 10^{-14}\) whereat the maximum also occurs in the boundary
cells.

\subsection{Periodic boundary conditions}
\label{sec:peri-bound-cond}

\begin{figure}[!htpb]
  \includegraphics[width=0.5\textwidth]{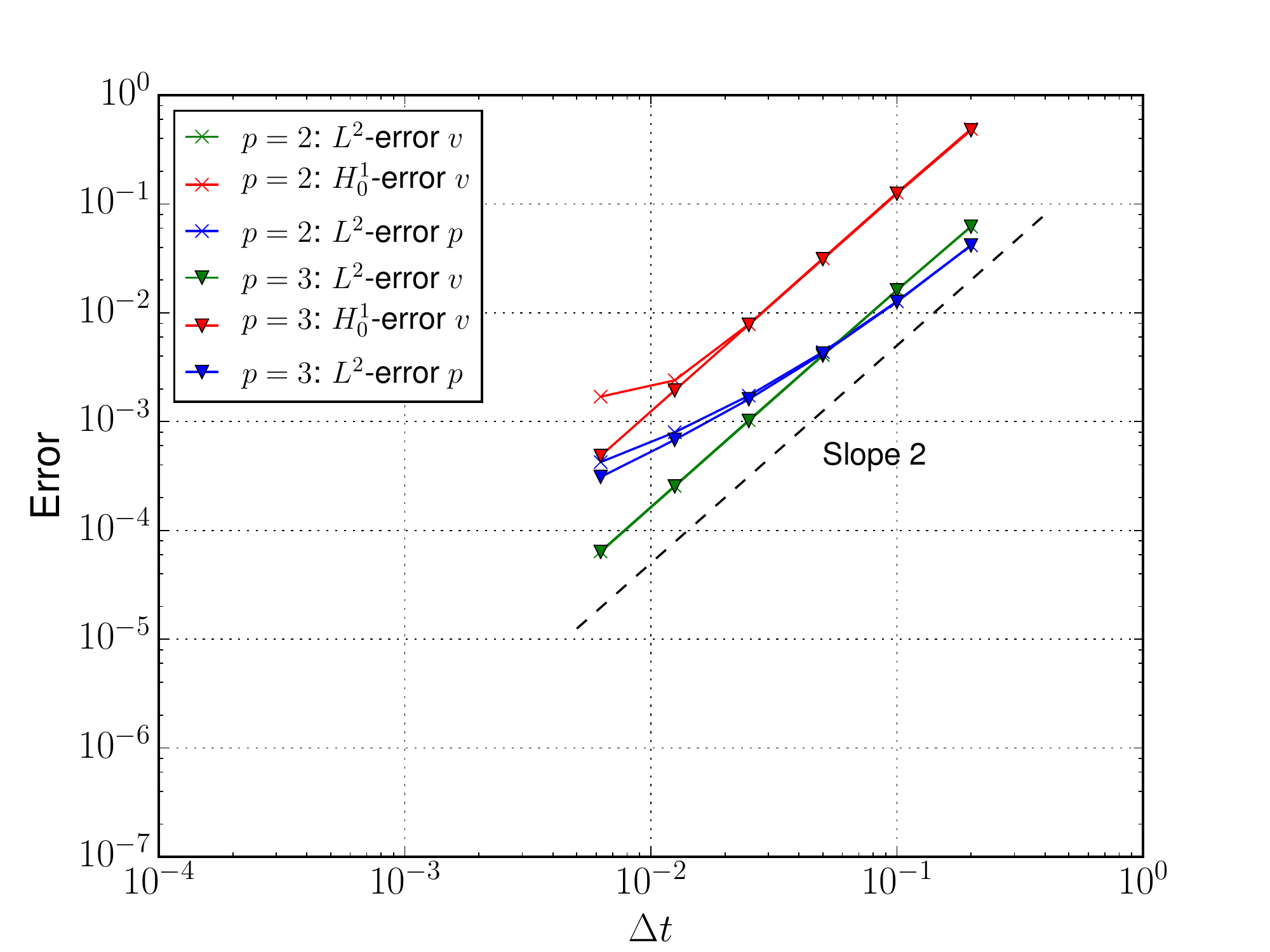}
  \includegraphics[width=0.5\textwidth]{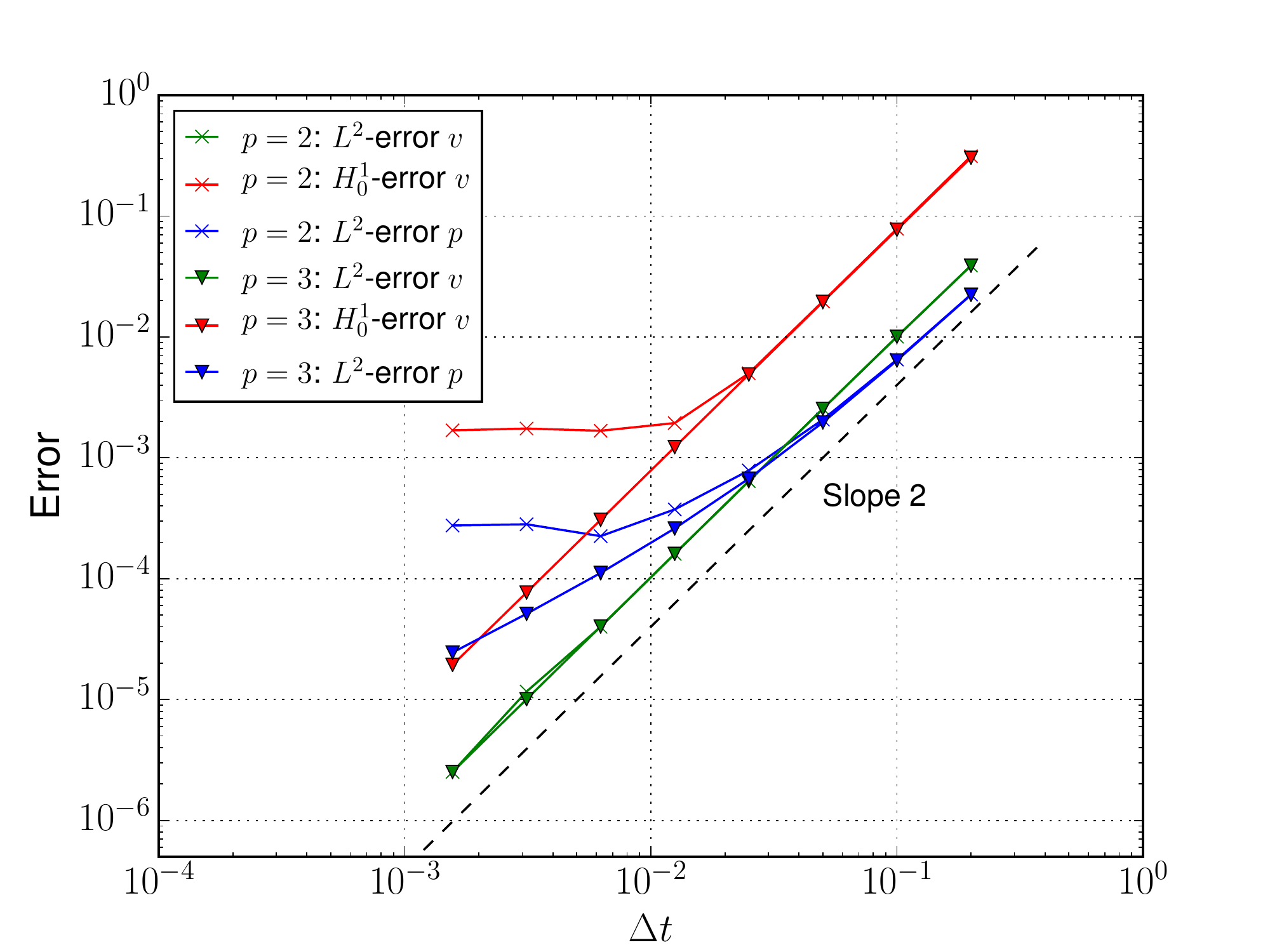}
  \caption{Errors and convergence rates at final time T=2 for the
    periodic boundary condition problem and spatial discretizations
    \THcube{2}{1}, \THcube{3}{2}, \THcube{4}{3}.
    Left part shows the IPCS.
    Right part shows the RIPCS.}
  \label{fig:ipcs_ripcs_periodic_rates}
\end{figure}

We continue with the configuration and test problem presented in
\ref{sec:local-mass-cons-1}. Figure \ref{fig:ipcs_ripcs_periodic_rates}
shows the error and convergence rates as a function of \(\Delta t\) for
the IPCS and RIPCS.
The green curves show the \(L^2\)-error for the velocity, the red
curves the \(H_0^1\)-error for the velocity and the blue curves the
\(L^2\)-error for the pressure obtained by the polynomial degrees
\(p=2,3\). The results for \(p=4\) are almost identical to \(p=3\),
therefore it has been omitted for the sake of clarity. For the IPCS the
curves grouped by the same color are almost identical meaning that the
splitting error is dominant in the measured range of \(\Delta
t\). Note however that for \(p=2\) in the RIPCS the spatial error
is already not negligible in this range and becomes all-dominant for
additionally smaller time steps taken. It can be seen on the right that
the \(H_0^1\)-error on the velocity and \(L^2\)-error on the pressure
flattens out whereas the errors from spatial order three continue
decreasing with the same rate. That puts in favour higher polynomial
degrees since the error on the same spatial mesh for moderate time step
sizes is minimized.

There is no rigorous error analysis of the projection methods for purely
periodic boundary conditions. But since in the periodic case no
artificial boundary conditions are imposed on the pressure, both the
standard and rotational formulation are expected to be fully second
order accurate. This is validated for the pressure-correction schemes in
figure \ref{fig:ipcs_ripcs_periodic_rates}.
The error of the RIPCS is slightly lower than the error of the IPCS, but
both schemes have the same convergence rate. It is close to
$\mathcal{O}(\Delta t^2)$ in the $L^2$-norm on the pressure while the
rates of the velocity in the $L^2$-norm and $H_0^1$-norm are perfectly
of second order.

The absence of artificial boundary conditions also implies that
the error on local mass conservation is distributed over the
interior on the domain. This was shown before in figure
\ref{pic:ripcs_order2_mass_conservation}.

\subsection{Beltrami flow}
\label{sec:beltrami-flow}

The Beltrami flow is one of the rare test problems where an exact fully
three-dimensional solution of the Navier-Stokes equations is derived. It
has its origin from \cite{EthierSteinmann1994} and has been later
studied by \cite{fenics:book}. The domain is \(\Omega = (-1,1)^3\) and
global Dirichlet boundary conditions are imposed by the exact solution
\begin{equation}
\begin{split}
  v_1(x,y,z,t) &= -a\,e^ {- d^2\,t }\,\left(e^{a\,x}\,\sin \left(d\,z+a\,y\right)+
 \cos \left(d\,y+a\,x\right)\,e^{a\,z}\right) \\
  v_2(x,y,z,t) &= -a\,e^ {- d^2\,t }\,\left(e^{a\,x}\,\cos \left(d\,z+a\,y\right)+e^{
 a\,y}\,\sin \left(a\,z+d\,x\right)\right) \\
  v_3(x,y,z,t) &= -a\,e^ {- d^2\,t }\,\left(e^{a\,y}\,\cos \left(a\,z+d\,x\right)+
 \sin \left(d\,y+a\,x\right)\,e^{a\,z}\right) \\
  p(x,y,z,t) &= -0.5\,a^2\,\rho\,e^ {- d^2\,t }\,(2\,\cos \left(d\,y+a\,x
 \right)\,e^{a\,\left(z+x\right)}\,\sin \left(d\,z+a\,y\right) \\
  &\quad +2\,e^{ a\,\left(y+x\right)}\,\sin \left(a\,z+d\,x\right)\,\cos \left(d\,z+a
 \,y\right) \\
  &\quad +2\,\sin \left(d\,y+a\,x\right)\,e^{a\,\left(z+y\right)}\,
 \cos \left(a\,z+d\,x\right)+e^{2\,a\,z}+e^{2\,a\,y}+e^{2\,a\,x} ) \; .
\end{split}
\label{eq:beltrami_3dproblem}
\end{equation}
The Beltrami flow has the property that the velocity and vorticity
vectors are aligned, namely \(d\; v - \nabla\times v = 0\)
The source term is given by \(f=0\), the density, viscosity are set to
\(\rho = \mu = 1\). The constants \(a\) and \(d\) may be chosen arbitrarily and have been
set to \(a=\pi/4\), \(d=\pi/2\) as in
\cite{EthierSteinmann1994}. Computations were done on a \(50\times
50\times 50\) cubic mesh.

\begin{figure}[!htpb]
  \centering
  \includegraphics[width=0.5\textwidth]{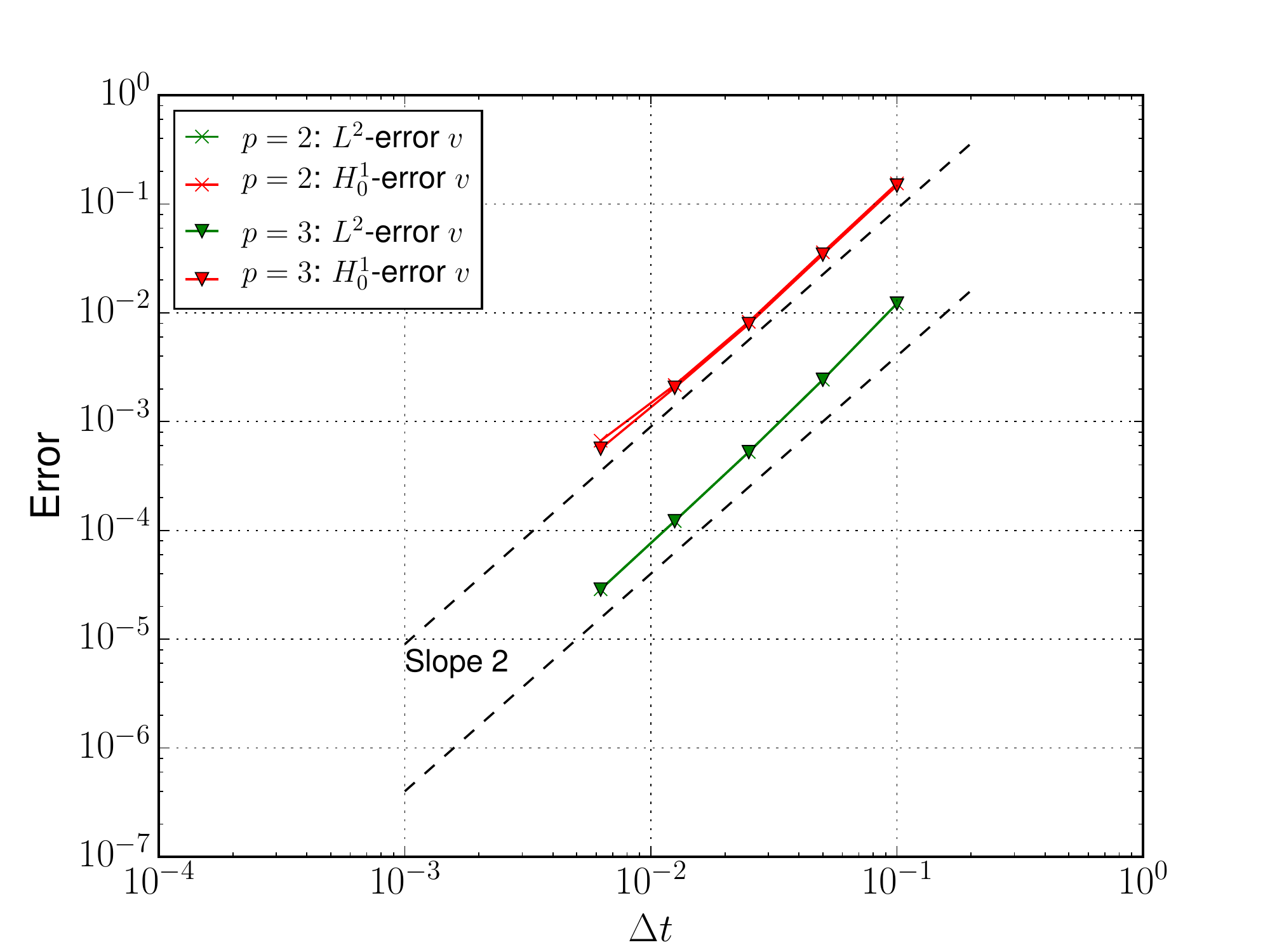}
  \caption{Errors and convergence rates at final time T=0.5 for the
    Beltrami test problem using the RIPCS with the spatial
    discretizations \THcube{2}{1}, \THcube{3}{2}.}
  \label{fig:ripcs_beltrami_rates}
\end{figure}

Figure \ref{fig:ripcs_beltrami_rates} shows the error and convergence
rates as a function of \(\Delta t\) for the RIPCS. The green curves show
the \(L^2\)-error for the velocity and the red curves the \(H_0^1\)-error
for the velocity obtained by the polynomial degrees \(p=2,3\). The
curves grouped by the same color are almost identical meaning that the
splitting error is dominant in the measured range of \(\Delta t\). It
can be concluded from the figure that error is fully second order
convergent in both norms.

\subsection{3D DNS of turbulent flows}
\label{sec:3d-dns-turbulent}

We test the applicability of the code to direct numerical simulations
(DNS) with two examples.

\subsubsection{3D Driven Cavity}
\label{sec:3d-driven-cavity}

The driven cavity flow is a commonly used benchmark problem due to its
simple geometry. We take the configuration based on
\cite{1873-7005-5-3-A03}. We shift the domain to \(\Omega =
(0,1)^3\) instead, set \(\rho = 1\) and \(\textup{Re} = 1/\mu =
10000\). The driving velocity at \(y=1\) has a continuous ramp profile
in time given by \(v = (\min(t,1),0,0)^T\). The simulation is computed
on a \(50\times 50\times 50\) grid using the \THcube{3}{2}
discretization with reconstruction in \(\RT_h^1\).
The temporal interval starts from \(t_0 = 0\) up to 50,
Alexander's second order strongly S-stable scheme is used with time step
size \(\Delta t = 0.005\).

\begin{figure}[!htpb]
  \includegraphics[width=0.5\textwidth]{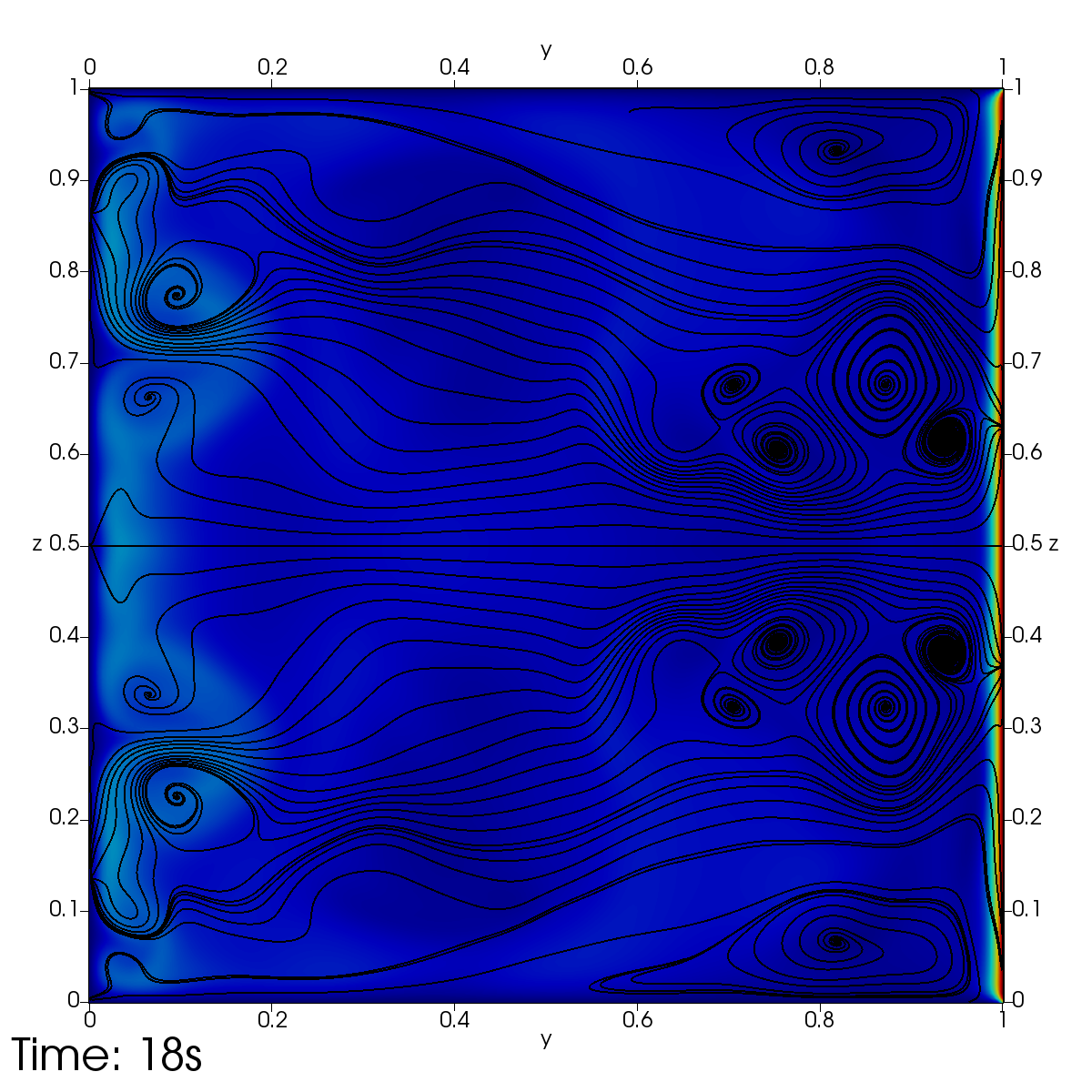}
  \includegraphics[width=0.5\textwidth]{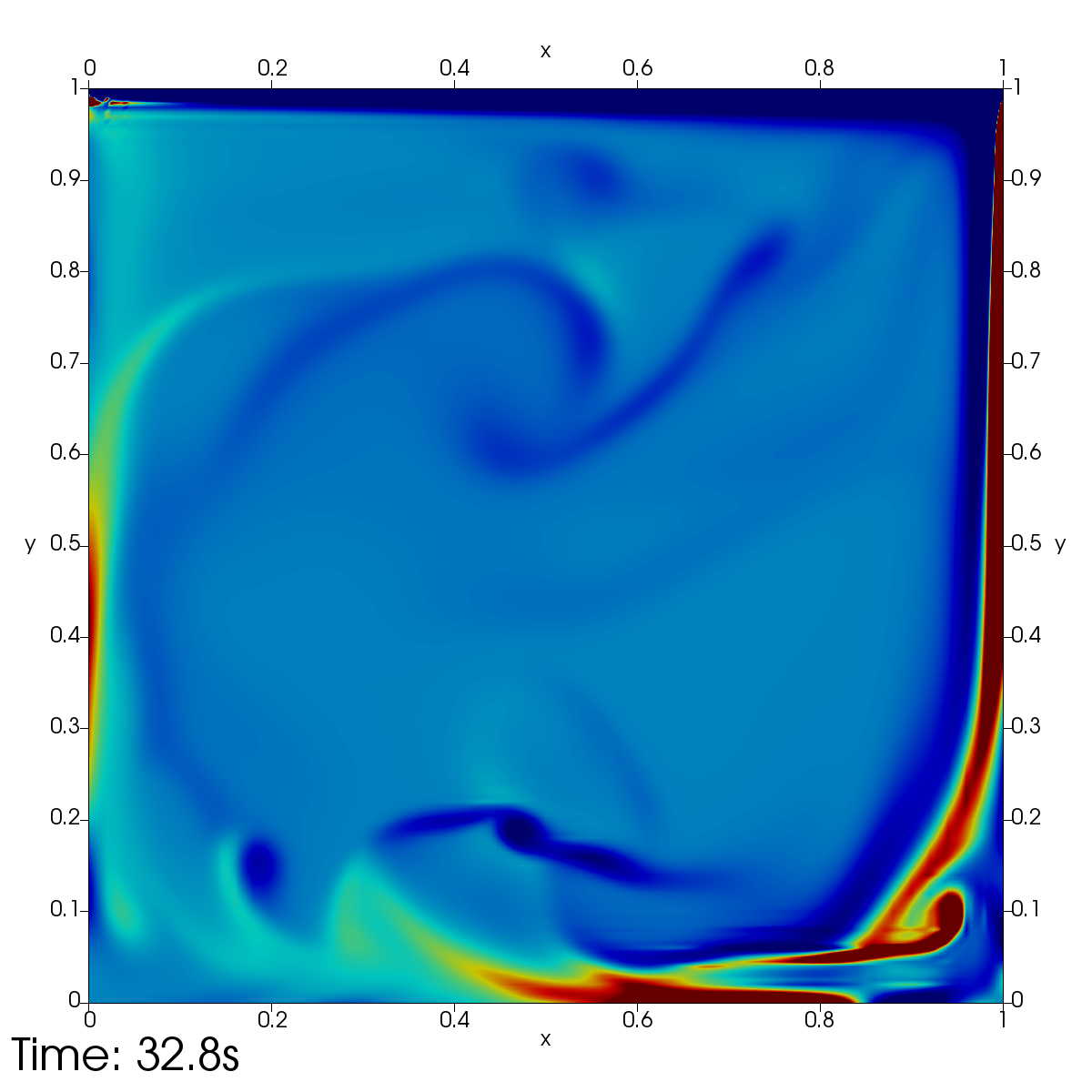}
  \caption{
    Time-snapshot of \(\lef v\ri_2\) and streamlines of \(v_\parallel\)
    on plane \(x = 0.5\) (left). Time-snapshot of z-component of the
    vorticity \(\nabla\times v\) on plane \(z=0.5\) (right).
  }
  \label{fig:driven_cavity_3d_re10000}
\end{figure}

The flow transitions to a chaotic behaviour at \(t\sim
12\). Figure \ref{fig:driven_cavity_3d_re10000} shows two snapshots of
the flow in the driven cavity. The left part of the figure shows \(\lef
v\ri_2\) together with streamlines of \(v_\parallel\) on the plane \(x =
0.5\). The characteristic corner vortices on this cut through the domain
are clearly visible as well as the appearance of the Taylor-G\"{o}rtler
vortices close to \(z=0.5\). The right part of the figure shows the
z-component of the vorticity \(\nabla\times v\) on the plane \(z=0.5\)
after the main initial vortex has decayed into several small eddies.

The figure demonstrates in general long time stability of the spatial
discretization \THcube{3}{2} and reconstruction in \(\RT_h^1\) with
\(10^4\) time steps taken.

\subsubsection{3D Taylor-Green vortex}
\label{sec:3d-taylor-green}

The Taylor-Green vortex has been studied in great detail in
\cite{TaylorGreen3D}, C3.5. It aims at testing the accuracy and the
performance of high-order methods in the DNS. The initial flow field is
given by
\begin{align}
  v_1(x,y,z,t) &= V_0 \sin\left(\frac{x}{L}\right)
  \cos\left(\frac{y}{L}\right) \cos\left(\frac{z}{L}\right) \notag \\
  v_2(x,y,z,t) &= -V_0 \cos\left(\frac{x}{L}\right)
  \sin\left(\frac{y}{L}\right) \cos\left(\frac{z}{L}\right)\notag \\
  v_3(x,y,z,t) &= 0 \notag \\
  p(x,y,z,t) &= p_0 + \frac{\rho_0 V_0^2}{16} \left(
    \cos\left(\frac{2x}{L}\right) + \cos\left(\frac{2y}{L}\right)\right)
  \left(\cos\left(\frac{2z}{L}\right) + 2\right)
  \label{eq:taylorgreen_vortex_3d_initial}
\end{align}
with periodic boundary conditions in all directions.  The general
extension of the domain is \(\Omega = (-\pi L, \pi L)^3\) and the
Reynolds number is given by \(\textup{Re} = \frac{\rho_0 V_0
  L}{\mu}\). As in the references \cite{fenics:book,TaylorGreen3D} we
set \(L=1, \; V_0 = 1, \; \rho_0 = 1, \; p_0 = 0, \; \textup{Re} =
1600\). In three dimensions the flow transitions to turbulence with
development of small scale structures. The simulation is computed on a
\(105\times 105\times 105\) grid using the \THcube{3}{2} discretization
with reconstruction in \(\RT_h^1\) which leads to about \(2\cdot 10^8\)
degrees of freedom. The temporal interval starts from \(t_0 = 0\) up to
20, Alexander's second order strongly S-stable scheme is used with time
step size \(\Delta t = 0.005\).

\begin{figure}[!htpb]
  \includegraphics[width=0.6\textwidth]{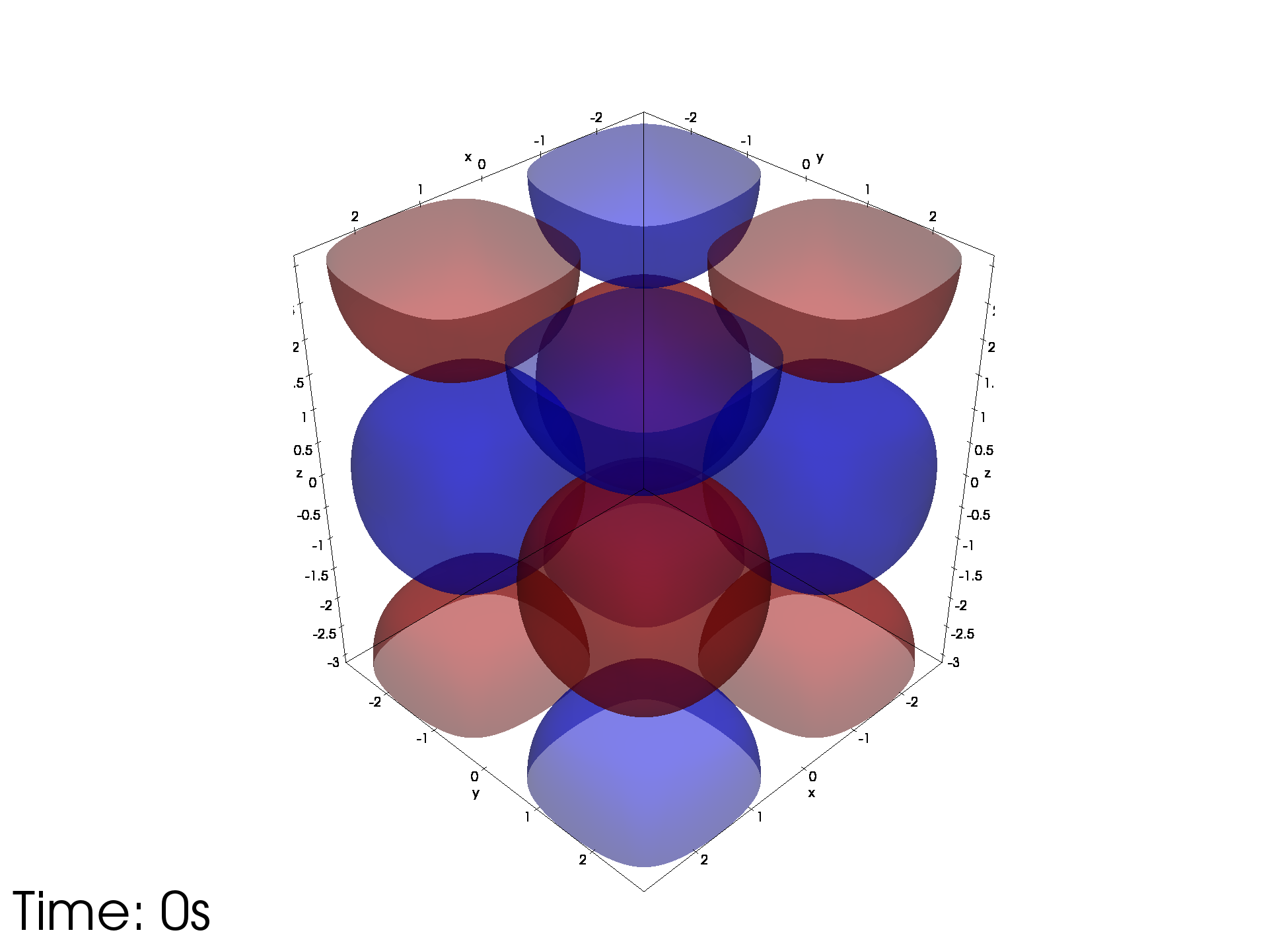}
  \hspace*{-3em}
  \includegraphics[width=0.6\textwidth]{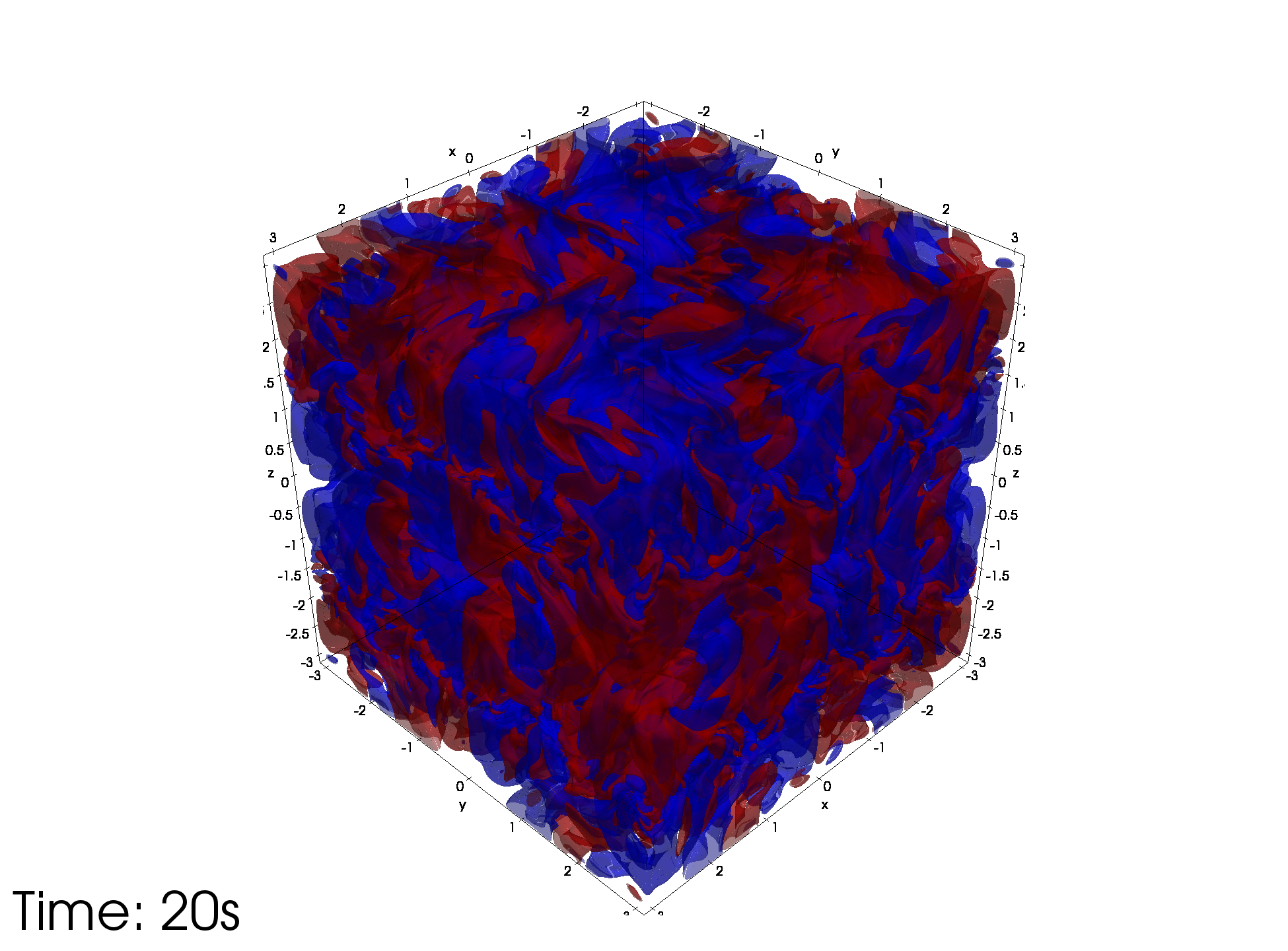}
  \caption{
    Vorticity visualization for the Taylor-Green vortex
    problem. The red contours show the levels of the vorticity
    z-component with 25 percent of the instantaneous peak. The blue
    contours show the corresponding negative level.
    }
  \label{fig:taylorgreen_vortex_3d_re1600}
\end{figure}

Figure \ref{fig:taylorgreen_vortex_3d_re1600} shows the contour surfaces
of the z-component of the vorticity for the values 0.5 in red and -0.5
in blue, respectively, at initial and final condition.

\begin{figure}[!htpb]
  \includegraphics[width=0.5\textwidth]{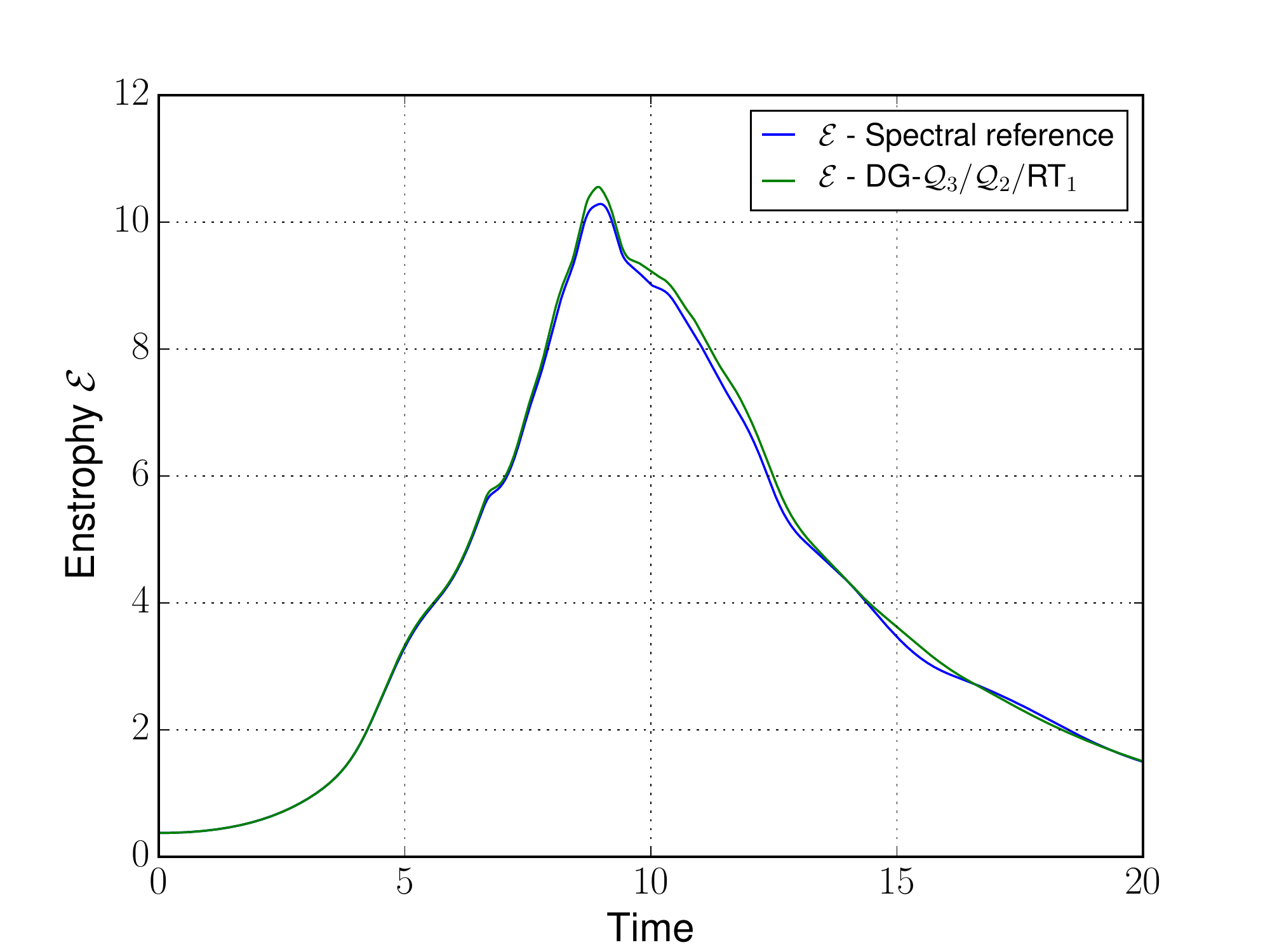}
  \includegraphics[width=0.5\textwidth]{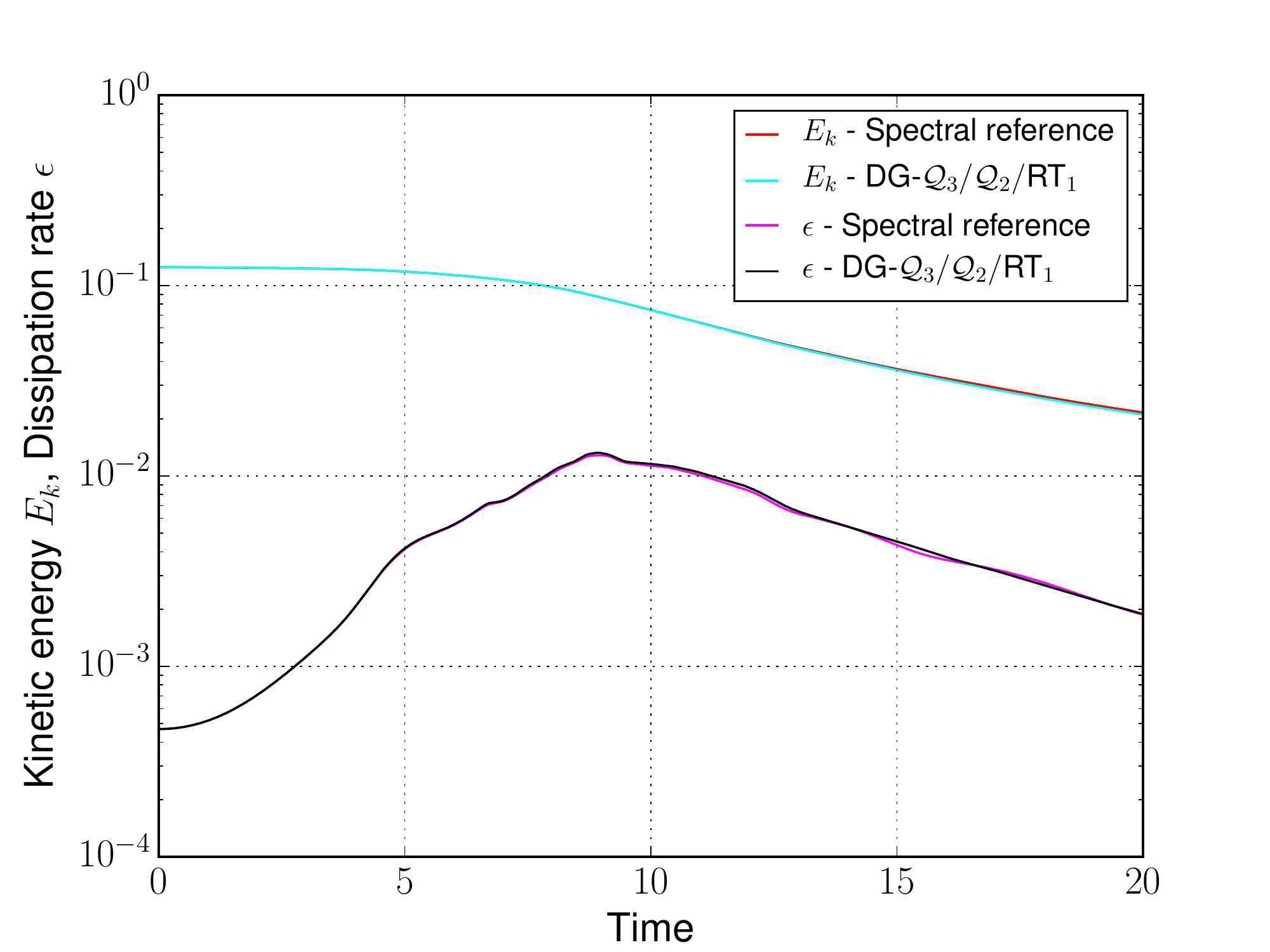}
  \caption{ Evolution of the kinetic energy, dissipation rate and
    enstrophy for the \THcube{3}{2} discretization with reconstruction
    in \(\RT_h^1\) and comparison with the reference values from the
    spectral code.}
  \label{fig:taylorgreen_vortex_3d_reference_values}
\end{figure}

We now compare our results to the reference values given for this
problem that contain the temporal evolution of
\begin{itemize}
\item the kinetic energy, \(E_k = \frac{1}{\rho \left|\Omega\right|} \frac12
  (v,v)_{0,\Omega} \),
\item the dissipation rate, \(\epsilon = \frac{\nu}{\left|\Omega\right|}
  (\nabla v, \nabla v)_{0,\Omega}\),
\item the enstrophy, \(\qquad \; \mathcal{E} = \frac{1}{\rho \left|\Omega\right|}
  \frac12 (\nabla\times v, \nabla\times v)_{0,\Omega} \; . \)
\end{itemize}
The reference solution was obtained with a dealiased pseudo-spectral code
run on a \(512^3\) grid, time integration was performed with a
low-storage three-step Runge-Kutta scheme and a time step of \(\Delta t
= 10^{-3}\). The comparison of these three reference quantities is
presented in figure
\ref{fig:taylorgreen_vortex_3d_reference_values}. Our enstrophy curve is
extremely close to the reference curve, furthermore the curves for the
kinetic energy and the dissipation rate are indistinguishable with
respect to the reference curves.

It can be concluded that the spatial discretization \THcube{3}{2} and
reconstruction in \(\RT_h^1\) together with the upwind scheme based on
the Vijayasundaram flux exhibits long time stability - as in the
three-dimensional driven cavity problem - and is good at capturing
turbulence accurately.

\section{Conclusion and Outlook}
\label{sec:concluding-remarks}

In this paper we considered splitting methods, in particular the incremental and rotational incremental pressure
correction schemes in combination with a high-order discontinuous Galerkin discretization in space.
The momentum equation is solved fully-implicitly with a matrix-free Newton method and
sum-factorized element computations to achieve reduced computational complexity at
high floating-point performance. The upwind discretization of the convective term in the momentum
equation uses a modified Vijayasundaram numerical flux function that takes into account
that the discrete velocity field is not in $H(\text{div})$.
In the Helmholtz projection step we employ an $H(\text{div})$
postprocessing of the velocity correction which ensures that the projected velocity satisfies the
discrete continuity equation. Numerical results confirm that RIPCS
is second-order convergent in time for Dirichlet and periodic boundary conditions and
has convergence order 1.5 for mixed boundary conditions. Three-dimensional computations with up
to $2\cdot 10^8$ degrees of freedom and about $10^4$ time steps show that the scheme
is stable for long time computations.

The $H(\text{div})$ postprocessing in $RT_k$ spaces is restricted to
parallelepiped elements due to the Piola transformation. It remains an open problem
how to extend the postprocessing  scheme to more general element transformations.
In a forthcoming publication we will focus on the performance characteristics and
scalability of the parallel implementation.

%

\section*{Acknowledgements}

The first author (M. P.) is supported by a Ph.D. stipend of the
Heidelberg Graduate School of Mathematical
and Computational Methods for the Sciences (GSC 220). The second
author (S. M.) is funded
within the DFG special program Software for Exascale Computing (SPP 1648) under contract number
Ba 1498/10-2.  Computing resources were provided by bwHPC supported by the state of Baden-Württemberg.
We would also like to thank Eike M\"{u}ller for supplying us with a
matrix-free block Jacobi and block SOR code for linear equations.

\bibliographystyle{elsarticle-num}
\bibliography{literature}

\end{document}

%% file: convergence-tables/ripcs_order2_periodic_rates.out.tex
\begin {tabular}{c|c|c|c}%
dt&$L^2$ error $v$&$H_0^1$ error $v$&$L^2$ error $p$\\\hline %
2.000e-01&3.89336e-02&3.12769e-01&2.24301e-02\\%
1.000e-01&1.00536e-02&7.86088e-02&6.49292e-03\\%
5.000e-02&2.54833e-03&1.96963e-02&2.07180e-03\\%
2.500e-02&6.40802e-04&5.00184e-03&7.85322e-04\\%
1.250e-02&1.60275e-04&1.93785e-03&3.74214e-04\\%
6.250e-03&3.99586e-05&1.67405e-03&2.24884e-04\\\hline %
\end {tabular}%

%% file: convergence-tables/ripcs_hdiv_order2_periodic_rates.out.tex
\begin {tabular}{c|c|c|c}%
dt&$L^2$ error $v$&$H_0^1$ error $v$&$L^2$ error $p$\\\hline %
2.000e-01&3.89026e-02&3.06282e-01&2.24215e-02\\%
1.000e-01&1.00444e-02&7.72283e-02&6.49243e-03\\%
5.000e-02&2.54548e-03&1.96643e-02&2.07260e-03\\%
2.500e-02&6.39597e-04&5.44801e-03&7.86167e-04\\%
1.250e-02&1.59573e-04&2.63178e-03&3.75384e-04\\%
6.250e-03&3.97373e-05&2.33037e-03&2.27288e-04\\\hline %
\end {tabular}%

%% file: convergence-tables/ripcs_rt1_order2_periodic_rates.out.tex
\begin {tabular}{c|c|c|c}%
dt&$L^2$ error $v$&$H_0^1$ error $v$&$L^2$ error $p$\\\hline %
2.000e-01&3.88824e-02&3.03453e-01&2.24163e-02\\%
1.000e-01&1.00387e-02&7.71449e-02&6.49254e-03\\%
5.000e-02&2.54417e-03&1.96504e-02&2.07356e-03\\%
2.500e-02&6.39432e-04&5.44560e-03&7.86989e-04\\%
1.250e-02&1.59699e-04&2.63182e-03&3.75972e-04\\%
6.250e-03&3.99134e-05&2.33056e-03&2.27698e-04\\\hline %
\end {tabular}%

%% file: convergence-tables/ripcs_order3_periodic_rates.out.tex
\begin {tabular}{c|c|c|c}%
dt&$L^2$ error $v$&$H_0^1$ error $v$&$L^2$ error $p$\\\hline %
2.000e-01&3.88929e-02&3.03675e-01&2.23165e-02\\%
1.000e-01&1.00411e-02&7.70569e-02&6.38400e-03\\%
5.000e-02&2.54558e-03&1.94734e-02&1.96199e-03\\%
2.500e-02&6.40610e-04&4.89254e-03&6.72911e-04\\%
1.250e-02&1.60670e-04&1.22611e-03&2.60387e-04\\%
6.250e-03&4.02320e-05&3.06909e-04&1.11552e-04\\\hline %
\end {tabular}%

%% file: convergence-tables/ripcs_hdiv_order3_periodic_rates.out.tex
\begin {tabular}{c|c|c|c}%
dt&$L^2$ error $v$&$H_0^1$ error $v$&$L^2$ error $p$\\\hline %
2.000e-01&3.88969e-02&3.03581e-01&2.23175e-02\\%
1.000e-01&1.00423e-02&7.70650e-02&6.38402e-03\\%
5.000e-02&2.54588e-03&1.94755e-02&1.96185e-03\\%
2.500e-02&6.40684e-04&4.89306e-03&6.72797e-04\\%
1.250e-02&1.60689e-04&1.22624e-03&2.60320e-04\\%
6.250e-03&4.02364e-05&3.06934e-04&1.11516e-04\\\hline %
\end {tabular}%

%% file: paper.bbl
\begin{thebibliography}{10}
\expandafter\ifx\csname url\endcsname\relax
  \def\url#1{\texttt{#1}}\fi
\expandafter\ifx\csname urlprefix\endcsname\relax\def\urlprefix{URL }\fi
\expandafter\ifx\csname href\endcsname\relax
  \def\href#1#2{#2} \def\path#1{#1}\fi

\bibitem{girault_riviere}
B.~Rivière, V.~Girault,
  \href{http://www.sciencedirect.com/science/article/pii/S0045782505002690}{Discontinuous
  finite element methods for incompressible flows on subdomains with
  non-matching interfaces}, Computer Methods in Applied Mechanics and
  Engineering 195~(25–28) (2006) 3274 -- 3292, discontinuous Galerkin
  Methods.
\newblock \href {http://dx.doi.org/http://dx.doi.org/10.1016/j.cma.2005.06.014}
  {\path{doi:http://dx.doi.org/10.1016/j.cma.2005.06.014}}.
\newline\urlprefix\url{http://www.sciencedirect.com/science/article/pii/S0045782505002690}

\bibitem{Girault_adiscontinuous}
V.~Girault, B.~Rivière, Mary, F.~Wheeler, A discontinuous galerkin method with
  nonoverlapping domain decomposition for the stokes and navier-stokes
  problems, Math. Comp (2004) 53--84.

\bibitem{Ferrer2011224}
E.~Ferrer, R.~Willden,
  \href{http://www.sciencedirect.com/science/article/pii/S0045793010002860}{A
  high order discontinuous galerkin finite element solver for the
  incompressible navier–stokes equations}, Computers \& Fluids 46~(1) (2011)
  224 -- 230, 10th \{ICFD\} Conference Series on Numerical Methods for Fluid
  Dynamics (ICFD 2010).
\newblock \href
  {http://dx.doi.org/http://dx.doi.org/10.1016/j.compfluid.2010.10.018}
  {\path{doi:http://dx.doi.org/10.1016/j.compfluid.2010.10.018}}.
\newline\urlprefix\url{http://www.sciencedirect.com/science/article/pii/S0045793010002860}

\bibitem{2016arXiv160701323K}
B.~{Krank}, N.~{Fehn}, W.~A. {Wall}, M.~{Kronbichler}, {A high-order
  semi-explicit discontinuous Galerkin solver for 3D incompressible flow with
  application to DNS and LES of turbulent channel flow}, ArXiv e-prints\href
  {http://arxiv.org/abs/1607.01323} {\path{arXiv:1607.01323}}.

\bibitem{sum-fact-dg}
S.~M\"{u}thing, M.~Piatkowski, P.~Bastian, High-performance implementation of
  matrix-free spectral discontinuous galerkin methods, In preparation.

\bibitem{CBO9780511574856A009}
E.~J. Dean, R.~Glowinski,
  \href{http://dx.doi.org/10.1017/CBO9780511574856.003}{On some finite element
  methods for the numerical simulation of incompressible viscous flow}, in:
  M.~D. Gunzburger, R.~A. Nicolaides (Eds.), Incompressible Computational Fluid
  Dynamics, Cambridge University Press, 1993, pp. 17--66, cambridge Books
  Online.
\newline\urlprefix\url{http://dx.doi.org/10.1017/CBO9780511574856.003}

\bibitem{Elman}
H.~C. Elman, D.~J. Silvester, A.~J. Wathen, Finite elements and fast iterative
  solvers: with applications in incompressible fluid dynamics, Oxford
  University Press, 2014, second edition.

\bibitem{ChorinProjectionNavierStokes}
A.~J. Chorin, \href{http://www.jstor.org/stable/2004575}{Numerical solution of
  the navier-stokes equations}, Mathematics of Computation 22~(104) (1968)
  745--762.
\newline\urlprefix\url{http://www.jstor.org/stable/2004575}

\bibitem{TemamProjectionNavierStokes}
R.~Témam, \href{http://dx.doi.org/10.1007/BF00247696}{Sur l'approximation de
  la solution des équations de navier-stokes par la méthode des pas
  fractionnaires (ii)}, Archive for Rational Mechanics and Analysis 33~(5)
  (1969) 377--385.
\newblock \href {http://dx.doi.org/10.1007/BF00247696}
  {\path{doi:10.1007/BF00247696}}.
\newline\urlprefix\url{http://dx.doi.org/10.1007/BF00247696}

\bibitem{Rannacher1992}
R.~Rannacher, \href{http://dx.doi.org/10.1007/BFb0090341}{On chorin's
  projection method for the incompressible navier-stokes equations}, Springer
  Berlin Heidelberg, Berlin, Heidelberg, 1992, pp. 167--183.
\newblock \href {http://dx.doi.org/10.1007/BFb0090341}
  {\path{doi:10.1007/BFb0090341}}.
\newline\urlprefix\url{http://dx.doi.org/10.1007/BFb0090341}

\bibitem{WeinanE1}
W.~E, J.-G. Liu, Projection method i: Convergence and numerical boundary
  layers, SIAM Journal on Numerical Analysis 32~(4) (1995) 1017--1057.
\newblock \href {http://dx.doi.org/10.1137/0732047}
  {\path{doi:10.1137/0732047}}.

\bibitem{WeinanE2}
W.~E, J.-G. Liu, Projection method ii: Godunov–ryabenki analysis, SIAM
  Journal on Numerical Analysis 33~(4) (1996) 1597--1621.
\newblock \href {http://dx.doi.org/10.1137/S003614299426450X}
  {\path{doi:10.1137/S003614299426450X}}.

\bibitem{KARNIADAKIS1991414}
G.~E. Karniadakis, M.~Israeli, S.~A. Orszag,
  \href{http://www.sciencedirect.com/science/article/pii/0021999191900078}{High-order
  splitting methods for the incompressible navier-stokes equations}, Journal of
  Computational Physics 97~(2) (1991) 414 -- 443.
\newblock \href
  {http://dx.doi.org/http://dx.doi.org/10.1016/0021-9991(91)90007-8}
  {\path{doi:http://dx.doi.org/10.1016/0021-9991(91)90007-8}}.
\newline\urlprefix\url{http://www.sciencedirect.com/science/article/pii/0021999191900078}

\bibitem{FLD:FLD373}
L.~J.~P. Timmermans, P.~D. Minev, F.~N. Van De~Vosse, An approximate projection
  scheme for incompressible flow using spectral elements, International Journal
  for Numerical Methods in Fluids 22~(7) (1996) 673--688.

\bibitem{GuermondShen}
J.~L. Guermond, J.~Shen, {A new class of truly consistent splitting schemes for
  incompressible flows}, J. Comput. Phys. 192 (2003) 262--276.

\bibitem{Guermond04onthe}
J.~L. Guermond, J.~Shen, On the error estimates of rotational
  pressure-correction projection methods, Math. Comp 73 (2004) 1719--1737.

\bibitem{Karniadakis2005}
G.~Karniadakis, S.~Sherwin, Spectral/hp Element Methods for Computational Fluid
  Dynamics, Oxford University Press, 2005.

\bibitem{Guermond06}
J.~L. Guermond, P.~Minev, J.~Shen, {An overview of projection methods for
  incompressible flows}, Computer Methods in Applied Mechanics and Engineering
  195 (2006) 6011--6045.

\bibitem{GODA197976}
K.~Goda,
  \href{http://www.sciencedirect.com/science/article/pii/0021999179900883}{A
  multistep technique with implicit difference schemes for calculating two- or
  three-dimensional cavity flows}, Journal of Computational Physics 30~(1)
  (1979) 76 -- 95.
\newblock \href
  {http://dx.doi.org/http://dx.doi.org/10.1016/0021-9991(79)90088-3}
  {\path{doi:http://dx.doi.org/10.1016/0021-9991(79)90088-3}}.
\newline\urlprefix\url{http://www.sciencedirect.com/science/article/pii/0021999179900883}

\bibitem{Steinmoeller2013480}
D.~Steinmoeller, M.~Stastna, K.~Lamb,
  \href{http://www.sciencedirect.com/science/article/pii/S0021999113004026}{A
  short note on the discontinuous galerkin discretization of the pressure
  projection operator in incompressible flow}, Journal of Computational Physics
  251 (2013) 480 -- 486.
\newblock \href {http://dx.doi.org/http://dx.doi.org/10.1016/j.jcp.2013.05.036}
  {\path{doi:http://dx.doi.org/10.1016/j.jcp.2013.05.036}}.
\newline\urlprefix\url{http://www.sciencedirect.com/science/article/pii/S0021999113004026}

\bibitem{Joshi2016120}
S.~M. Joshi, P.~J. Diamessis, D.~T. Steinmoeller, M.~Stastna, G.~N. Thomsen,
  \href{http://www.sciencedirect.com/science/article/pii/S0045793016301311}{A
  post-processing technique for stabilizing the discontinuous pressure
  projection operator in marginally-resolved incompressible inviscid flow},
  Computers \& Fluids 139 (2016) 120 -- 129, 13th \{USNCCM\} International
  Symposium of High-Order Methods for Computational Fluid Dynamics - A special
  issue dedicated to the 60th birthday of Professor David Kopriva.
\newblock \href
  {http://dx.doi.org/http://dx.doi.org/10.1016/j.compfluid.2016.04.021}
  {\path{doi:http://dx.doi.org/10.1016/j.compfluid.2016.04.021}}.
\newline\urlprefix\url{http://www.sciencedirect.com/science/article/pii/S0045793016301311}

\bibitem{BastianRiviere}
P.~Bastian, B.~Rivière,
  \href{http://dx.doi.org/10.1002/fld.562}{Superconvergence and h(div)
  projection for discontinuous galerkin methods}, International Journal for
  Numerical Methods in Fluids 42~(10) (2003) 1043--1057.
\newblock \href {http://dx.doi.org/10.1002/fld.562}
  {\path{doi:10.1002/fld.562}}.
\newline\urlprefix\url{http://dx.doi.org/10.1002/fld.562}

\bibitem{ErnHdiv2007}
A.~Ern, S.~Nicaise, M.~Vohralík,
  \href{http://www.sciencedirect.com/science/article/pii/S1631073X07004360}{An
  accurate {H}(div) flux reconstruction for discontinuous {G}alerkin
  approximations of elliptic problems}, Comptes Rendus Mathematique 345~(12)
  (2007) 709 -- 712.
\newblock \href
  {http://dx.doi.org/http://dx.doi.org/10.1016/j.crma.2007.10.036}
  {\path{doi:http://dx.doi.org/10.1016/j.crma.2007.10.036}}.
\newline\urlprefix\url{http://www.sciencedirect.com/science/article/pii/S1631073X07004360}

\bibitem{doi:10.1137/040604418}
J.~L. Guermond, P.~Minev, J.~Shen,
  \href{http://dx.doi.org/10.1137/040604418}{Error analysis of
  pressure-correction schemes for the time-dependent stokes equations with open
  boundary conditions}, SIAM Journal on Numerical Analysis 43~(1) (2005)
  239--258.
\newblock \href {http://arxiv.org/abs/http://dx.doi.org/10.1137/040604418}
  {\path{arXiv:http://dx.doi.org/10.1137/040604418}}, \href
  {http://dx.doi.org/10.1137/040604418} {\path{doi:10.1137/040604418}}.
\newline\urlprefix\url{http://dx.doi.org/10.1137/040604418}

\bibitem{GiraultRaviartBook}
V.~Girault, P.~Raviart, Finite Element Methods for the Navier-Stokes Equations,
  Springer, 1986.

\bibitem{TemamBook}
R.~Témam, Navier-Stokes Equations. Theory and numerical analysis, North
  Holland, Amsterdam, 1987.

\bibitem{HoustonHartmann2008}
P.~Houston, R.~Hartmann, An optimal order interior penalty discontinuous
  {G}alerkin discretization of the compressible {N}avier-{S}tokes equations, J.
  Comp. Phys. 227 (2008) 9670--9685.

\bibitem{amg4dg}
P.~Bastian, M.~Blatt, R.~Scheichl, Algebraic multigrid for discontinuous
  {G}alerkin discretizations, Numer. Linear Algebra Appl. 19~(2) (2012)
  367--388.
\newblock \href {http://dx.doi.org/10.1002/nla.1816}
  {\path{doi:10.1002/nla.1816}}.

\bibitem{feistauer_book}
M.~Feistauer, J.~Felcman, I.~Straskraba, {Mathematical and Computational
  Methods for Compressible Flow}, Clarendon Press, 2003.

\bibitem{Dole2004727}
V.~Dolej\v{s}\'{i}, M.~Feistauer,
  \href{http://www.sciencedirect.com/science/article/pii/S0021999104000609}{A
  semi-implicit discontinuous galerkin finite element method for the numerical
  solution of inviscid compressible flow}, Journal of Computational Physics
  198~(2) (2004) 727 -- 746.
\newblock \href {http://dx.doi.org/http://dx.doi.org/10.1016/j.jcp.2004.01.023}
  {\path{doi:http://dx.doi.org/10.1016/j.jcp.2004.01.023}}.
\newline\urlprefix\url{http://www.sciencedirect.com/science/article/pii/S0021999104000609}

\bibitem{Evans}
L.~Evans, Partial Differential Equations, 2nd Edition, American Mathematical
  Society, 2010.

\bibitem{book:941322}
J.-L. G.~a. Alexandre~Ern, Theory and Practice of Finite Elements, 1st Edition,
  Applied Mathematical Sciences 159, Springer-Verlag New York, 2004.

\bibitem{Schweizer2013}
B.~Schweizer, Partielle Differentialgleichungen, Eine anwendungsorientierte
  Einf\"{u}hrung, Springer-Verlag, 2013.

\bibitem{Bhatia:2013:HDS:2498747.2498999}
H.~Bhatia, G.~Norgard, V.~Pascucci, P.-T. Bremer,
  \href{http://dx.doi.org/10.1109/TVCG.2012.316}{The helmholtz-hodge
  decomposition\&\#x2014;a survey}, IEEE Transactions on Visualization and
  Computer Graphics 19~(8) (2013) 1386--1404.
\newblock \href {http://dx.doi.org/10.1109/TVCG.2012.316}
  {\path{doi:10.1109/TVCG.2012.316}}.
\newline\urlprefix\url{http://dx.doi.org/10.1109/TVCG.2012.316}

\bibitem{Brezzi91}
F.~Brezzi, M.~Fortin, Mixed and Hybrid Finite Element Methods, Springer-Verlag,
  1991.

\bibitem{Alexander:1977:DIR}
R.~Alexander, Diagonally implicit {Runge--Kutta} methods for stiff
  {O}.{D}.{E}.'s, SIAM Journal on Numerical Analysis 14~(6) (1977) 1006--1021.

\bibitem{dune08:1}
P.~Bastian, M.~Blatt, A.~Dedner, C.~Engwer, R.~Kl{\"o}fkorn, M.~Ohlberger,
  O.~Sander, A generic grid interface for parallel and adaptive scientific
  computing. {P}art {I}: {A}bstract framework, Computing 82~(2--3) (2008)
  103--119.
\newblock \href {http://dx.doi.org/10.1007/s00607-008-0003-x}
  {\path{doi:10.1007/s00607-008-0003-x}}.

\bibitem{Bastian2016}
P.~Bastian, C.~Engwer, J.~Fahlke, M.~Geveler, D.~Göddeke, O.~Iliev,
  O.~Ippisch, R.~Milk, J.~Mohring, S.~Müthing, M.~Ohlberger, D.~Ribbrock,
  S.~Turek,
  \href{http://dx.doi.org/10.1007/978-3-319-40528-5_1}{{Hardware-Based
  Efficiency Advances in the EXA-DUNE Project}}, Springer International
  Publishing, Cham, 2016, p. 3–23.
\newblock \href {http://dx.doi.org/10.1007/978-3-319-40528-5_1}
  {\path{doi:10.1007/978-3-319-40528-5_1}}.
\newline\urlprefix\url{http://dx.doi.org/10.1007/978-3-319-40528-5_1}

\bibitem{Kronbichler2012}
M.~Kronbichler, K.~Kormann, A generic interface for parallel cell-based finite
  element operator application, Computers \& Fluids 63 (2012) 135--147.

\bibitem{10.2307/96892}
G.~I. Taylor, A.~E. Green, \href{http://www.jstor.org/stable/96892}{Mechanism
  of the production of small eddies from large ones}, Proceedings of the Royal
  Society of London. Series A, Mathematical and Physical Sciences 158~(895)
  (1937) 499--521.
\newline\urlprefix\url{http://www.jstor.org/stable/96892}

\bibitem{5982828}
P.~Rabenold, E.~Balaras, Parallel adaptive mesh refinement for the
  incompressible navier-stokes equations.

\bibitem{EthierSteinmann1994}
C.~R. Ethier, D.~A. Steinmann, Exact fully {3D} {N}avier--{S}tokes solution for
  benchmarking, Internat. J. Numer. Methods Fluids 19 (1994) 369 -- 375.

\bibitem{fenics:book}
A.~Logg, K.-A. Mardal, G.~N. Wells (Eds.),
  \href{http://dx.doi.org/10.1007/978-3-642-23099-8}{Automated Solution of
  Differential Equations by the Finite Element Method}, Vol.~84 of Lecture
  Notes in Computational Science and Engineering, Springer, 2012.
\newblock \href {http://dx.doi.org/10.1007/978-3-642-23099-8}
  {\path{doi:10.1007/978-3-642-23099-8}}.
\newline\urlprefix\url{http://dx.doi.org/10.1007/978-3-642-23099-8}

\bibitem{1873-7005-5-3-A03}
R.~Iwatsu, K.~Ishii, T.~Kawamura, K.~Kuwahara, J.~M. Hyun,
  \href{http://stacks.iop.org/1873-7005/5/i=3/a=A03}{{Numerical simulation of
  three-dimensional flow structure in a driven cavity}}, Fluid Dynamics
  Research 5~(3) (1989) 173.
\newline\urlprefix\url{http://stacks.iop.org/1873-7005/5/i=3/a=A03}

\bibitem{TaylorGreen3D}
\href{http://www.dlr.de/as/hiocfd}{2nd international workshop on high-order cfd
  methods, cologne, germany} (May 2013).
\newline\urlprefix\url{http://www.dlr.de/as/hiocfd}

\end{thebibliography}
